\documentclass[preprint,11pt]{elsarticle}
\usepackage{mathrsfs}
\usepackage{amsfonts}
\usepackage{amssymb}
\usepackage{color}
\usepackage{algorithm}
\usepackage{multirow}
\usepackage{xcolor}
\usepackage{setspace}
\usepackage{epstopdf}

\usepackage{subfigure}

\usepackage{booktabs}
\usepackage[style=base]{caption}
\biboptions{numbers,sort,compress}
\usepackage{latexsym,mathdots,array,extarrows}
\providecommand{\U}[1]{\protect\rule{.1in}{.1in}}
\allowdisplaybreaks[4]
\journal{Applied Numerical Mathematics}

\usepackage{algpseudocode,float}

\usepackage{caption}
\usepackage{graphicx}
\makeatletter\def\CT{\def\@captype{figure}}\makeatother



\numberwithin{equation}{section}
\newenvironment{AMS}{\noindent{\bf Mathematics Subject Classification:}}

\newtheorem{theorem}{Theorem}
\newtheorem{example}{Example}
\newtheorem{lemma}[theorem]{Lemma}

\newtheorem{definition}{Definition}
\newtheorem{remark}{Remark}

\newtheorem{proposition}{Proposition}

\newenvironment{proof}{\noindent{\bf Proof:}}{\hfill\fbox{}\vspace*{1mm}}

\begin{document}	
	\begin{frontmatter}
		\title{A parallel-in-time preconditioner for Crank-Nicolson discretization of a parabolic optimal control problem}
		\author[add1]{Xue-Lei Lin\corref{cor1}}\ead{linxuelei@hit.edu.cn}  
		\author[add2]{Shu-Lin Wu}\ead{wushulin84@hotmail.com}  		    	
		\cortext[cor1]{The Corresponding Author.}
		\address[add1]{School of Science, Harbin Institute of Technology,	Shenzhen 518055, China}
		\address[add2]{School of Mathematics and Statistics, Northeast Normal University, Changchun 130024, China}
		
		\begin{abstract}
			In this paper, a fast solver is studied for saddle point system arising from a second-order Crank-Nicolson discretization of an initial-valued parabolic PDE constrained  optimal control problem, which is indefinite and ill-conditioned. Different from the saddle point system arising from the first-order Euler discretization, the saddle point system arising from Crank-Nicolson discretization has a dense and non-symmetric Schur complement, which brings challenges to fast solver designing. To remedy this, a novel symmetrization technique is applied to the saddle point system so that the new Schur complement is symmetric definite and the well-known matching-Schur-complement (MSC) preconditioner is applicable to the new Schur complement. Nevertheless, the new Schur complement is still a dense matrix and the inversion of the corresponding MSC preconditioner is  not parallel-in-time (PinT) and thus time consuming. For this concern, a modified MSC preconditioner for the new Schur complement system. Our new preconditioner can be implemented in a fast and PinT way via a temporal diagonalization technique. Theoretically,  the eigenvalues of the preconditioned matrix by our new preconditioner are proven to be lower and upper bounded by positive constants independent of matrix size and the regularization parameter.  With such spectrum, the preconditioned conjugate gradient (PCG) solver for the Schur complement system is proven to have a convergence rate independent of matrix size and regularization parameter. To the best of my knowledge, it is the first time to have an iterative solver with problem-independent convergence rate for the saddle point system arising from Crank-Nicolson discretization of the optimal control problem.
			Numerical results are reported to show that the performance of the proposed preconditioner.
		\end{abstract}
		
	\end{frontmatter}
	
	\section{Introduction}\label{introduction}
	Consider the 
	optimal control problem:
	\begin{equation}\label{optcontrlprob}
		\min\limits_{y,u} L(y,u):=\frac{1}{2}||y-g||_{L^2(\Omega\times(0,T))}^2+\frac{\gamma}{2}||u||_{L^2(\Omega\times(0,T))}^2,
	\end{equation}
	subject to a parabolic equation with initial- and boundary-value  conditions 
	\begin{equation*}
		\begin{cases}
			y_{t}- Ly=f+u, \quad {\rm in~~}\Omega\times(0,T),\\
			y=0,\quad{\rm on~~}\partial\Omega\times(0,T),\\
			y(\cdot,0)=y_0,\quad \Omega,
		\end{cases}
	\end{equation*}
	where $g$, $f$, $y_0$  are all given functions and $\gamma>0$ is a given regularization parameter. The goal of the optimization process is to find the state $y$ as close as possible to the desired state $g$ using the control $u$. The space operator $ L$ considered here is the diffusion operator $ L=\nabla(a(x)\nabla\cdot)$ with $a(x)\geq 0$, but other self-adjoint operators can be also included.

	The optimal control problem \eqref{optcontrlprob}, constrained by a state equation in the form of a parabolic equation arises in many applications, see, e.g., \cite{wang2004simplified,troltzsch2010optimal,kunisch2015time}.  Since the analytical solutions of optimal control problems are typically unavailable, the problems are most often solved numerically. Numerical methods for optimal control problems mainly consists of two classes, i.e., discretize-then-optimize method and optimize-then-discretize method, see, e.g., \cite{rao2009survey,troltzsch2010optimal,polak1973historical,miele1975recent,von1992direct,leitmann1981calculus,donald2016optimal} and the references therein. Following the optimize-then-discretize method, the optimal control problem \eqref{optcontrlprob} is firstly converted into the following KKT system by applying the first-order optimality condition (see, e.g., \cite{abbeloos2011nested,biegler2007real} for more details):
	\begin{equation*}
		\left[
		\begin{array}[c]{ccc}
			I&0& L_2\\
			0&\gamma I&- I\\
			L_1&- I&0
		\end{array}
		\right]\left[\begin{array}[c]{c}
			y\\
			u\\
			p
		\end{array}\right]=\left[\begin{array}[c]{c}
			g\\
			0\\
			f
		\end{array}\right],
	\end{equation*}
	where $p$ is the Lagrange multiplier; $ L_1$ and $ L_2$ are defined by
	\begin{align*}
		& L_1y:=(\partial_t - L)y,~ {\rm with}~y({\bf x},0)=y_0({\bf x}),\\
		& L_2p:=(-\partial_t- L)p,~ {\rm with}~p({\bf x},T)=0.
	\end{align*}
	Here, $y$ evolves forward from $t=0$ to $t=T$ and $p$ evolves backward from $t=T$ to $t=0$. By eliminating the control variable $u$, the KKT system is reduced to
	\begin{equation}\label{reducedkkt}
		\left[
		\begin{array}[c]{cc}
			I& L_2\\
			L_1&-\gamma^{-1} I
		\end{array}
		\right]\left[\begin{array}[c]{c}
			y\\
			p
		\end{array}\right]=\left[\begin{array}[c]{c}
			g\\
			f
		\end{array}\right].
	\end{equation}
	This work is focused on studying fast numerical solution for \eqref{reducedkkt}.
	
	In order to solve \eqref{reducedkkt}, numerical discretization is then applied to converting the problem into an approximate problem formulated in finite dimension space. For space discretization of \eqref{reducedkkt}, the central difference scheme is adopted, which leads to a symmetric positive semi-definite  matrix $ L_{h}$ he Laplacian  operator $ L$  is approximated by a symmetric positive semi-definite  matrix $ L_{h}$ as discretization of $-L$. To match the second-order accuracy of the spatial discretization, the second-order-accuracy Crank-Nicolson discretization  \cite{apel2012crank} is employed for approximating the temporal derivative. The resulting  discrete KKT system is 
	\begin{equation}\label{origdiscrdckktsys}
		\left[\begin{array}[c]{cc}
			\frac{\tau}{2} B_2\otimes I_J&B_1^T\otimes I_J+\frac{\tau}{2}B_2^T\otimes L_{h}\\
			B_1\otimes I_J+\frac{\tau}{2}B_2\otimes L_{h}&-\frac{\tau}{2\gamma}B_2^{\rm T}\otimes I_J
		\end{array}\right]\left[\begin{array}[c]{c}
			y_{\tau,h}\\
			p_{\tau,h}
		\end{array}\right]=\left[\begin{array}[c]{c}
			g_{\tau,h}\\
			f_{\tau,h}
		\end{array}\right],
	\end{equation}
	where $I_k$ denotes the $k\times k$ identity matrix;  $J$ is the number of spatial grid points; $\tau=N/T$ is the time step size with $N$ being the number of temporal grid points;  The two matrices $B_1$ and $B_2$ are 
	\begin{equation}\label{b1b2def}
		B_1=\left[\begin{array}[c]{cccc}
			1&&&\\
			-1&1&&\\
			&\ddots&\ddots& \\
			&&-1&1
		\end{array}\right],\quad B_2=\left[\begin{array}[c]{cccc}
			1&&&\\
			1&1&&\\
			&\ddots&\ddots& \\
			&&1&1
		\end{array}\right]\in\mathbb{R}^{N \times N}.
	\end{equation}
	The symbol `$\otimes$' denotes the Kronecker product. This paper focus on studying a fast solver for  the linear system \eqref{origdiscrdckktsys}. 
	
	Typically, the coefficient matrix of \eqref{origdiscrdckktsys} is of large size and sparse and \eqref{origdiscrdckktsys} must be solved iteratively. Nevertheless,  \eqref{origdiscrdckktsys} is a type of saddle point problem \cite{benzi2005numerical}, which is  ill-conditioned and indefinite \footnote{The indefiniteness here means the symmetric part of the matrix is indefinite}. Because of this, classical iterative solvers converges slowly for saddle point problems. To deal with the issue, preconditioning techniques are proposed to improve the spectral distribution of the saddle point problems, see, e.g., \cite{rees2010preconditioning,cao2016simplified,wang2009optimization,notay2014new,pearson2012regularization,wathen1995convergence,elman2014finite,wu2020diagonalization,liu2020parameter,ipsen2001note,de2005block,axelsson2003preconditioning,cao2016simplified,franceschini2019block,quirynen2020presas,wu2020parallel,sogn2019schur} and the references therein. Among them, the preconditioning techniques based on Schur complement significantly reduce the  iteration number of Krylov subspace solvers for the preconditioned  systems, see, e.g., \cite{de2005block,ipsen2001note,pearson2012regularization,loghin2004analysis,murphy2000note}. The implementation of the Schur complement based preconditioning techniques involves inversion of the Schur complement. Nevertheless, unlike the symmetric  saddle point system arising from first-order accuracy temporal discretization scheme, the saddle point system \eqref{origdiscrdckktsys} arising from the second-order Crank-Nicolson is nonsymmetric, which leads to a  non-symmetric and dense  Schur complement. Such complicated Schur complement  lacks of suitable algorithms for fast inversion. To address this issue, a novel block-diagonal scaling process is applied to symmetrizing \eqref{origdiscrdckktsys}, which is described as follows.
	
	Clearly, \eqref{origdiscrdckktsys} is equivalent to
	\begin{align}
		&\left[\begin{array}[c]{cc}
			\frac{\tau}{2} B_2\otimes I_J&B_1^T\otimes I_J+\frac{\tau}{2}B_2^T\otimes L_{h}\\
			B_1\otimes I_J+\frac{\tau}{2}B_2\otimes L_{h}&-\frac{\tau}{2\gamma}B_2^{\rm T}\otimes I_J
		\end{array}\right]W^{-1}\left[\begin{array}[c]{c}
			\tilde{y}_{\tau,h}	\\
			\tilde{p}_{\tau,h}
		\end{array}\right]=\left[\begin{array}[c]{c}
			g_{\tau,h}\\
			f_{\tau,h}
		\end{array}\right],\label{scaledkktsys1}\\
		&\left[\begin{array}[c]{c}
			y_{\tau,h}	\\
			p_{\tau,h}
		\end{array}\right]=W^{-1}\left[\begin{array}[c]{c}
			\tilde{y}_{\tau,h}\\
			\tilde{p}_{\tau,h}
		\end{array}\right],\label{scallingsystem}
	\end{align}
	where  $W={\rm blockdiag}(B_2\otimes I_J,B_2^{\rm T}\otimes I_J)$. It is straightforward to verify that
	\begin{align}
		&B_2^{-1}=\left[
		\begin{array}
			[c]{cccc}
			s_{0} &       &  &\\
			s_{1}& s_{0} &  &\\
			\vdots&\ddots &\ddots&\\
			s_{N-1}& \ldots & s_{1} & s_{0}
		\end{array}
		\right],\quad s_k=(-1)^{k},~k=0,1,...,N-1,\label{b2inversionform}\\	
		&B_1B_2^{-1}=B_2^{-1}B_1=\left[
		\begin{array}
			[c]{cccc}
			q_{0} &       &  &\\
			q_{1}& q_{0} &  &\\
			\vdots&\ddots &\ddots&\\
			q_{N-1}& \ldots & q_{1} & q_{0}
		\end{array}
		\right],\quad q_0=1,\quad q_k=2(-1)^{k},~ k\geq 1.\label{b1b2invcommuteq}
	\end{align}
	\eqref{b2inversionform} and \eqref{b1b2invcommuteq} show that $B_2^{-1}$, $B_2^{-1}B_1$ are both Toeplitz matrices and that $B_1$ and $B_2^{-1}$ are commutable.
	Since   $W^{-1}={\rm blockdiag}(B_2^{-1}\otimes I_J,B_2^{-\rm T}\otimes I_J)$ and $B_2^{-1}$ has a simple expression as shown in \eqref{b2inversionform},  $p_{\tau,h}$ and $y_{\tau,h}$ can be fast computed from \eqref{scallingsystem} within $\mathcal{O}(JN\log N)$ flops by means of fast Fourier transforms (FFTs) once $\tilde{y}_{\tau,h}$ and $\tilde{p}_{\tau,h}$ are computed.
	
	Denote 
	$$
	B=B_2^{-1}B_1.
	$$
	Then, the commutativity between $B_1$ and $B_2^{-1}$ imply that
	\begin{align*}
		&\left[\begin{array}[c]{cc}
			\frac{\tau}{2} B_2\otimes I_J&B_1^T\otimes I_J+\frac{\tau}{2}B_2^T\otimes L_{h}\\
			B_1\otimes I_J+\frac{\tau}{2}B_2\otimes L_{h}&-\frac{\tau}{2\gamma}B_2^{\rm T}\otimes I_J
		\end{array}\right]W^{-1}\\
		&=\left[\begin{array}[c]{cc}
			\frac{\tau}{2} I_N\otimes I_J&	B^{\rm T}\otimes I_J+\frac{\tau}{2}I_N\otimes L_{h}\\
			B\otimes I_J+\frac{\tau}{2}I_N\otimes L_{h}&	-\frac{\tau}{2\gamma}I_N\otimes I_J
		\end{array}\right],
	\end{align*}
	where $I_N\in\mathbb{R}^{N\times N}$ denotes the $N\times N$ identity matrix.
	Now one can see that \eqref{origdiscrdckktsys} is equivalent to the following symmetric linear system 
	\begin{equation}\label{scaledkktsys}
		\left[\begin{array}[c]{cc}
			\frac{\tau}{2} I_N\otimes I_J&	B^{\rm T}\otimes I_J+\frac{\tau}{2}I_N\otimes L_{h}\\
			B\otimes I_J+\frac{\tau}{2}I_N\otimes L_{h}&	-\frac{\tau}{2\gamma}I_N\otimes I_J
		\end{array}\right]\left[\begin{array}[c]{c}
			\tilde{y}_{\tau,h}\\
			\tilde{p}_{\tau,h}
		\end{array}\right]=\left[\begin{array}[c]{c}
			g_{\tau,h}\\
			f_{\tau,h}
		\end{array}\right]. 
	\end{equation}
	Unlike the Schur complement of \eqref{origdiscrdckktsys}, the new  Schur complement  of \eqref{scaledkktsys} is
	\begin{equation}\label{schucompofscaledsys}
		-\frac{\tau}{2\gamma}I_N\otimes I_J-(B\otimes I_J+\frac{\tau}{2}I_N\otimes L_{h})\left(\frac{2}{\tau}I_N\otimes I_J\right)(B\otimes I_J+\frac{\tau}{2}I_N\otimes L_{h})^{\rm T},
	\end{equation}
	which is symmetric positive-definite (SPD),up to a multiplicative constant $-1$. With the Schur complement decomposition formula for inversion of $2\times 2$ block matrix, one can  see that the dominant complexity for solving the symmetrized saddle point system \eqref{scaledkktsys} is devoted to solving a linear system  with the Schur complement as coefficient matrix. Such linear system is called  Schur complement linear system throughout this paper. In \cite{pearson2012new}, an MSC preconditioner was proposed for the Schur complement system. According to the theory developed in \cite{pearson2012new},  the eigenvalues of the preconditioned Schur complement matrix are lower bounded by $0.5$ and upper bounded by 1 provided that  $B+B^{\rm T}$ is  positive semi-definite.   $B+B^{\rm T}$ is indeed positive semi-definite by the lateral analysis given in  Section \ref{spectralanalysis}. That means the PCG solver with the MSC preconditioner for solving the Schur complement system  has a linear convergence rate independent of matrix size and regularization parameter. Nevertheless, since $B$ is a dense triangular matrix, the inversion of the MSC preconditioner is  not PinT. In order to further reduce the computational time,  we in this work focus on developing a modified MSC preconditioner so that the inversion of the preconditioner is PinT and thus is faster to implement.
	
	The studies on PinT methods have gained increasing attention during last few decades, which aim to solve time-dependent problems in PinT manners; see, e.g., \cite{gander201550,lions2001661,falgout2014parallel,gander2007analysis,wu2018toward,wu2019acceleration,kwok2019schwarz,gander2016analysis,hackbusch1985parabolic,horton1995space,vandewalle2013parallel,liu2022well,liu2022rom} and the references therein. Paradiag method is a class of PinT methods designed based on diagonalization of (approximated) discrete temporal operators appearing in time-dependent problems; see, e.g., \cite{gander2016direct,mcdonald2018preconditioning,mcdonald2017preconditioning,wathen2019note,gander2019convergence,lin2021all,gander2023new,gander2020paradiag,wu2020parallel,wu2020diagonalization,hon2022optimal,hon2022sine,lin2021parallel,lin2018separable,gu2020parallel}. Inspired by the Paradiag method based on block $\alpha$-circulant approximation developed in the literature (see, e.g., \cite{lin2018separable,lin2021all,wu2020parallel,liu2020fast}),  a modification to the MSC preconditioner is proposed in this work to improve the efficiency of inverting the MSC preconditioner.
	The modified MSC preconditioner is based on introducing block $\alpha$-circulant approximation to the factors of the MSC preconditioner. The advantage of such modification is that the inversions of the modified factors are fast and PinT, thanks to the fast block diagonalizability of the block $\alpha$-circulant matrices (See the discussion in Section \ref{precdefsection}). Such fast invertibility guarantees that each iteration of PCG solver with the modified MSC preconditioner can be fast implemented.  Moreover, utilizing approximation property of the modified MSC to the MSC preconditioner and the clustering property of the preconditioned matrix by MSC preconditioner,  the eigenvalues of the preconditioned matrix by the modified MSC preconditioner are proven to be lower and upper bounded by positive constants independent of  matrix size and regularization parameter provided that $\alpha$ is sufficiently small. Such spectral distribution guarantees that the PCG solver with the modified MSC preconditioner for solving the Schur complement system  has a convergence rate independent of matrix-size and parameter. That means the PCG solver  with the modified MSC preconditioner for solving the Schur complement system converges within a fixed iteration number  no matter how the parameters of the problem change. The uniformly bounded iteration number and the low operation cost at iteration guarantee the low operation cost of the proposed PCG method in total.

	The contribution of this paper is three fold. Firstly, a novel symmetrization technique is proposed for the saddle point system \eqref{origdiscrdckktsys} arising from the Crank-Nicolson discretization of the KKT system \eqref{reducedkkt} so that the MSC preconditioner and its corresponding theory are applicable to preconditioning the new Schur complement. Secondly, a modification  to the MSC preconditioner is proposed so that the inversion of the new preconditioner is fast  and PinT. Thirdly, with our proposed preconditioner, the PCG solver for the preconditioned system is proven to  have a convergence rate independent of matrix size and regularization parameter. 
	
	The rest of this paper is structured as follows. In Section 2, the Schur complement system, the modified MSC preconditioner and its  implementation are presented. In Section 3, the spectrum of the preconditioned matrix and optimal convergence of PCG solver for the preconditioned system are analyzed. In Section 4, numerical results are presented to show the performance of the proposed preconditioning technique. Finally, concluding remarks are given in Section 5.
	
	\section{The Schur Complement of \eqref{scaledkktsys} and the Proposed PinT Preconditioner}\label{precdefsection}
	In this section, the Schur complement system corresponding to \eqref{scaledkktsys} is presented and the PinT preconditioner is introduced for the Schur complement system.

	It is straightforward to verify that the unknown of \eqref{scaledkktsys} can be expressed follows
	\begin{align}
		&\left[\begin{array}[c]{c}
			\tilde{y}_{\tau,h}\\
			\tilde{p}_{\tau,h}
		\end{array}\right]=\left[\begin{array}[c]{cc}
			I&\frac{G^{\rm T}}{-\tau} \\
			&I
		\end{array}\right]\left[\begin{array}[c]{cc}
			\frac{2}{\tau}  I_{N}\otimes I_J&\\
			&-2\gamma(\tau I_{N}\otimes I_J+\eta GG^{\rm T})^{-1}
		\end{array}\right]\left[\begin{array}[c]{cc}
			I&\\
			\frac{G}{-\tau} &I
		\end{array}\right]\left[\begin{array}[c]{c}
			g_{\tau,h}\\
			f_{\tau,h}
		\end{array}\right]\label{kktinversform},
	\end{align}
	where $I$ denotes the $(NJ)\times  (NJ)$ identity matrix.
	\begin{equation}\label{gmatdef}
		\eta=\frac{\gamma}{\tau}, ~G=2B\otimes I_J+\tau I_N\otimes L_{h}.
	\end{equation}
	From \eqref{kktinversform}, we see that the most heavy computational burden for solving the unknown is to compute the matrix vector product $(\tau I_{N}\otimes I_J+\eta GG^{\rm T})^{-1}b$ for some given vector $b$. In other words, the dominant computational burden for solving the unknown is to solve the following Schur complement system
	\begin{equation}\label{origschurcmpsys}
		{\bf K}v=b,
	\end{equation}
	where $b\in\mathbb{R}^{NJ\times 1}$ is some given vector; 
	\begin{equation}\label{kmatdef}
		{\bf K}=\tau  I_{NJ}+\eta  GG^{\rm T}.
	\end{equation}
	A routine calculation yields  
	\begin{equation}\label{kmatexpan}
		{\bf K}=(\tau I_{N}+4\eta BB^{\rm T})\otimes I_J+2\eta \tau(B+B^{\rm T})\otimes L_{h}+\tau^2\eta I_N\otimes L_{h}^2.
	\end{equation}
	
	To construct the preconditioner for ${\bf K}$, we firstly introduce an intermediate approximation ${\bf P}$ to ${\bf K}$.
	Denote
	\begin{equation}\label{Pmeddef}
		{\bf P}:={\bf R}{\bf R}^{\rm T},
	\end{equation}
	where
	\begin{equation}\label{rmatdef}
		{\bf R}=(\sqrt{\tau}I_N+2\sqrt{\eta}B)\otimes I_J+\tau\sqrt{\eta} I_N\otimes L_{h}.
	\end{equation}
	The matrix ${\bf P}$ is the so-called MSC preconditioner; see \cite{pearson2012new}.
	Note that the inversion of ${\bf P}$ requires the inversion of ${\bf R}$ and ${\bf R}^{\rm T}$. However, ${\bf R}$ (${\bf R}^{\rm T}$) has block lower (upper) triangular structure. The inversion of block lower (upper)  triangular matrix  involves block forward (backward)  substitution algorithm, which is not PinT. To remedy this, a modified MSC preconditioner ${\bf P}_{\alpha}$ is introduced as follows.
	\begin{equation}\label{pepsdef}
		{\bf P}_{\alpha}={\bf R}_{\alpha}{\bf R}_{\alpha}^{\rm T},
	\end{equation}
	where
	\begin{align}
		&{\bf R}_{\alpha}=(\sqrt{\tau}I_N+2\sqrt{\eta}B_{\alpha})\otimes I_J+\tau\sqrt{\eta}I_N\otimes L_{h},\label{repsdef}\\
		&B_{\alpha}=B+\alpha\tilde{B},\quad \tilde{B}=\left[\begin{array}[c]{cccc}
			0&q_{N-1}&\ldots&q_1\\
			&\ddots&\ddots&\vdots\\
			&&0&q_{N-1}\\
			&&&0
		\end{array}\right],\quad \alpha>0. \label{bepsdef}
	\end{align}
	Here, we recall that the above numbers $q_i$'s are defined in \eqref{b1b2invcommuteq}. We will prove the invertibility of ${\bf R}_{\alpha}$ and thus the invertibility of ${\bf P}_{\alpha}$ for some properly chosen $\alpha$ in the next section. From the definition of ${\bf P}_{\alpha}$, it is clear that ${\bf P}_{\alpha}$ is a symmetric positive definite matrix. Instead of solving \eqref{origschurcmpsys} directly, we apply PCG solver with  ${\bf P}_{\alpha}$ as preconditioner to solving the following preconditioned system
	\begin{equation}\label{precschurcmpsys}
		{\bf P}_{\alpha}^{-1}{\bf K}v={\bf P}_{\alpha}^{-1}b.
	\end{equation}
	
	In each PCG iteration, we need to compute some matrix vector products ${\bf P}_{\alpha}^{-1}{\bf K}w$ for some given vector $w\in\mathbb{R}^{NJ\times 1}$. One can compute the product ${\bf K}w$ first and then compute the product ${\bf P}_{\alpha}^{-1}({\bf K}w)$ once ${\bf K}w$ is computed. As ${\bf K}$ is sparse matrix, the computation of the product ${\bf K}w$ requires $\mathcal{O}(NJ)$ flops. Hence, we only need to consider the computation of ${\bf P}_{\alpha
	}^{-1}{\bf w}$ for a given vector ${\bf w}$. From ${\bf P}_{\alpha}^{-1}={\bf R}_{\alpha}^{-1}{\bf R}_{\alpha}^{-\rm T}$, we know that ${\bf P}_{\alpha
	}^{-1}{\bf w}$ can be computed by ${\bf P}_{\alpha
	}^{-1}{\bf w}={\bf R}_{\alpha}^{-1}({\bf R}_{\alpha}^{-\rm T}{\bf w})$. Hence, to fast implement the matrix-vector product of ${\bf P}_{\alpha}^{-1}$, it suffices to fast implement the matrix-vector products corresponding to ${\bf R}_{\alpha}^{-1}$ and ${\bf R}_{\alpha}^{-\rm T}$. We firstly discuss a fast PinT implementation for computing the product between ${\bf R}_{\alpha}^{-1}$ and an arbitrarily given vector ${\bf r}$.
	
	Actually, $B_{\alpha}$ is diagonalizable by means of Fourier transform (see, \cite{bini2005numerical}). The diagonalization formula is shown as follows.
	\begin{equation}\label{bepsdiagform}
		B_{\alpha}=D_{\alpha}^{-1}F_N\Lambda_{\alpha}F_N^{*} D_{\alpha},
	\end{equation}
	where 
	\begin{align}
		&F_N=\frac{1}{\sqrt{N}}\left[\theta_{N}^{(i-1)(j-1)}\right]_{i,j=1}^{N},\quad \theta_{N}=\exp\left(\frac{2\pi{\bf i}}{N}\right),\quad {\bf i}=\sqrt{-1},\label{fmatdef}\\	
		&D_{\alpha}={\rm diag}\left(\alpha^{\frac{i-1}{N}}\right)_{i=1}^{N},~ \Lambda_{\alpha}={\rm diag}(\lambda_{i-1}^{(\alpha)})_{i=1}^{N},~ \lambda_{k}^{(\alpha)}=\sum\limits_{j=0}^{N-1}q_j\alpha^{\frac{j}{N}}\theta_{N}^{-kj},~ k=0,1,...,N-1.\notag
	\end{align}
	Here, $F_N$ is called a Fourier transform matrix which is unitary. $F_N^{*}$ ($F_N$) multiplied with a vector is equivalent to Fourier transform (inverse Fourier transform) of the vector up to a scaling constant. From the definitions of $\{\lambda_k^{(\alpha)}\}_{k=0}^{N-1}$'s, it is clear that the $N$ numbers $\{\lambda_k^{(\alpha)}\}_{k=0}^{N-1}$ can be computed within one fast Fourier transform (FFT), which requires only $\mathcal{O}(N\log N)$ flops.
	
	By the fact that $(D_{\alpha}^{-1}F_N)(F_N^{*} D_{\alpha})=I_N$ and \eqref{bepsdiagform}, we know that 
	\begin{equation}\label{repsbdiagform}
		{\bf R}_{\alpha}=[(D_{\alpha}^{-1}F_N)\otimes I_J][(\sqrt{\tau}I_N+2\sqrt{\eta}\Lambda_{\alpha})\otimes I_J+\tau\sqrt{\eta} I_N\otimes L_{h}][(F_N^{*}D_{\alpha})\otimes I_J].
	\end{equation}
	Then, for any given vector ${\bf r}\in\mathbb{R}^{NJ\times 1}$, the matrix-vector product ${\bf x}:={\bf R}_{\alpha}^{-1}{\bf r}$ can be computed by solving the linear system ${\bf R}_{\alpha}{\bf x}={\bf r}$ in the following PinT pattern
	\begin{equation}\label{threestepimplemt}
		\begin{cases}
			\tilde{\bf r}=[(F_N^{*}D_{\alpha})\otimes I_J]{\bf r}, &\text{step-(a)},\\
			\left((\sqrt{\tau}+ 2\sqrt{\eta}\lambda_{k-1}^{(\alpha)}){I}_x + \tau\sqrt{\eta} L_{h}\right)\tilde{x}_k=\tilde{r}_k, ~k=1,2,\dots, N, &\text{step-(b)},\\
			{\bf x}=[(D_{\alpha}^{-1}F_N)\otimes I_J]\tilde{\bf x},&\text{step-(c)}, 
		\end{cases}
	\end{equation}
	where $\tilde{\bf r}=(\tilde{r}_1^{\rm T},\tilde{r}_2^{\rm T},...,\tilde{r}_N^{\rm T})^{\rm T}$, $\tilde{\bf x}=(\tilde{x}_1^{\rm T},\tilde{x}_2^{\rm T},...,\tilde{x}_N^{\rm T})^{\rm T}$. 
	
	By applying FFTs along temporal dimension, we see that the computation of Step-(a) and Step-(c) of \eqref{threestepimplemt} requires only $\mathcal{O}(JN\log N)$ flops. The most heavy computation burden is to solve the $N$ many linear systems in Step-(b) of \eqref{threestepimplemt}. But interestingly, these systems are completely independent and therefore the computation is PinT.
	
	From \eqref{repsbdiagform}, we know that ${\bf R}_{\alpha}^{\rm T}$ is also block diagonalizable, i.e.,
	\begin{equation}\label{repstransbdiagform}
		{\bf R}_{\alpha}^{\rm T}=[(D_{\alpha}F_N)\otimes I_J][(\sqrt{\tau}I_N+2\sqrt{\eta}\bar{\Lambda}_{\alpha})\otimes I_J+\tau\sqrt{\eta} I_N\otimes L_{h}][(F_N^{*}D_{\alpha}^{-1})\otimes I_J].
	\end{equation}
	Hence, the matrix-vector product ${\bf R}_{\alpha}^{-\rm T}{\bf r}$ for a given vector ${\bf r}$ can be computed in a PinT manner similar to \eqref{threestepimplemt}.
	
	In summary, the whole  process for computing the matrix-vector product ${\bf P}_{\alpha}^{-1}{\bf w}$ with a given vector ${\bf w}\in\mathbb{R}^{NJ\times 1}$ is PinT, which is why  ${\bf P}_{\alpha}$ is called  a PinT preconditioner.
	\section{Spectra of the Preconditioned matrix ${\bf P}_{\alpha}^{-1}{\bf K}$ and Optimal Convergence of the PCG Solver for the Preconditioned System \eqref{precschurcmpsys}}\label{spectralanalysis}
	In this section, for proper choice of $\alpha$,  the eigenvalues of the preconditioned matrix are proven to be lower and upper bounded by constants in independent of $N$, $J$ and $\gamma$ . The uniform boundedness of the spectrum implies an optimal convergence rate of the PCG solver for the preconditioned system. 
	
	In what follows, the spectrum of ${\bf P}^{-1}{\bf K}$ and ${\bf P}_{\alpha}^{-1}{\bf P}$ are estimated  first and then the spectrum of the preconditioned matrix ${\bf P}_{\alpha}^{-1}{\bf K}$ is estimated based on the estimations of  spectrum of ${\bf P}^{-1}{\bf K}$ and ${\bf P}_{\alpha}^{-1}{\bf P}$.
	
	\begin{proposition}\label{pinvertbilityprop}
		${\bf P}$ is invertible.
	\end{proposition}
	\begin{proof}
		Since ${\bf P}={\bf R}{\bf R}^{\rm T}$, it suffices to show  the invertibility of ${\bf R}$. From \eqref{rmatdef}, we see that ${\bf R}$ is a block triangular matrix whose invertibility is determined by the invertibility of diagonal blocks. Notice that all diagonal blocks of ${\bf R}$ are identical to $(\sqrt{\tau}+2\sqrt{\eta})I_J+\tau\sqrt{\eta} L_{h}$ that is clearly positive definite and thus invertible. Hence, ${\bf R}$ is invertible. The proof is complete.
	\end{proof}

	For any real symmetric matrices ${\bf H}_1,{\bf H}_2\in\mathbb{R}^{m\times m}$, denote ${\bf H}_2 \succ ({\rm or} \succeq) \ {\bf H}_1$
	if ${\bf H}_2-{\bf H}_1$ is  real symmetric positive definite (or real symmetric positive semi-definite).
	Also, ${\bf H}_1 \prec ({\rm or} \preceq) \ {\bf H}_2$ has the same meaning as that of ${\bf H}_2 \succ ({\rm or} \succeq) \ {\bf H}_1$.
	
	Let $O$ denote zero matrix with proper size.

	For any $O\preceq {\bf H}\in\mathbb{R}^{m\times m}$ and $p\geq 0$, denote $${\bf H}^{p}:={\bf U}^{\rm T}{\rm diag}(d_1^{p},d_2^{p},...,d_m^{p}){\bf U},$$
	where ${\bf U}^{\rm T}{\rm diag}(d_1,d_2,...,d_m){\bf U}$ is orthogonal diagonalization of ${\bf H}$. In particular, if ${\bf H}\succ O$, then we write $({\bf H}^{-1})^{p}$ as ${\bf H}^{-p}$ for $p> 0$.
	
	Denote 
	\begin{equation}\label{dhatdef}
		\hat{D}={\rm diag}((-1)^{k})_{k=1}^{N}.
	\end{equation}
	Let $e_i$ denote the $i$th column of the $N\times N$ identity matrix. Denote $${\bf 1}_N=(1;1;\cdots;1)\in\mathbb{R}^{N\times 1}.$$ Then, it is straightforward to verify that
	\begin{equation}\label{bpbtransrank1form}
		B+B^{\rm T}=2\hat{D}{\bf 1}_N{\bf 1}_N^{\rm T}\hat{D}=2\hat{D}{\bf 1}_N{\bf 1}_N^{\rm T}\hat{D}^{-1},
	\end{equation}
	where the last equality is due to $\hat{D}=\hat{D}^{-1}$.

	For a real symmetric matrix $H$, let $\lambda_{\min}(H)$ and $\lambda_{\max}(H)$ denote the minimal and maximal eigenvalues of $H$, respectively. For a real square matrix ${\bf C}$, define its symmetric part $\mathcal{H}({\bf C})$ as
	\begin{equation*}
		\mathcal{H}({\bf C}):=\frac{{\bf C}+{\bf C}^{\rm T}}{2}.
	\end{equation*} 
	
	Since ${\bf P}_{\alpha}$ is used as a preconditioner, the invertibility of ${\bf P}_{\alpha}$ is a necessary property. To this end, we firstly show the invertibility of  ${\bf P}_{\alpha}$ in what follows.
	
	\begin{definition}\label{circmatdef}
		A Toeplitz matrix $C=[c_{i-j}]_{i,j=1}^{m}\in\mathbb{R}^{m\times m}$ is called a circulant matrix if and only if $c_{i}=c_{i-m}$ for $i=1,2,...,m-1$.
	\end{definition}

	Recall the Fourier transform matrix $F_{N}$ defined in \eqref{fmatdef}. Actually, any $N\times N$ circulant matrix is diagonalizable by $F_{N}$. We state this interesting fact in the following proposition.
	\begin{proposition}\textnormal{(see, e.g., \cite{ng2004iterative})}\label{circdiagprop}
		A circulant matrix $C\in\mathbb{C}^{N\times N}$ is diagonalizable by the $F_N$, i.e.,
		\begin{equation*}
			C=F_N{\rm diag}(\sqrt{N}F_N^{*}C(:,1))F_N^{*},
		\end{equation*}
		where $C(:,1)$ denotes the first column of $C$.
	\end{proposition}

	Recall the definition of ${\bf R}_{\alpha}$ given in \eqref{repsdef}, we see that
	\begin{equation}\label{hrepsexpan}
		\mathcal{H}({\bf R}_{\alpha})=(\sqrt{\tau}I_N+2\sqrt{\eta}\mathcal{H}(B_{\alpha}))\otimes I_J+\tau\sqrt{\eta}I_N\otimes L_{h}.
	\end{equation}
	From the definition of $B_{\alpha}$ given in \eqref{bepsdef}, it is straightforward to verify that
	\begin{equation}
		\mathcal{H}(B_{\alpha})=\hat{D}C_{\alpha}\hat{D}^{-1}=\hat{D}C_{\alpha}\hat{D},
	\end{equation}
	where $\hat{D}$ is a diagonal matrix defined in \eqref{dhatdef}, $C_{\alpha}$ is a circulant matrix with $$[1+\alpha(-1)^{N}]{\bf 1}_N-\alpha(-1)^Ne_1\in\mathbb{R}^{N\times 1},$$ as its first column. $e_1\in\mathbb{R}^{N\times 1}$ denotes the first column of the $N\times N$ identity matrix.
	
	Denote $\tilde{F}_N=\hat{D}F_N$. It is clear that $\tilde{F}_N$ is a unitary matrix. By Proposition \ref{circdiagprop}, $C_{\alpha}$ is diagonalizable by $F_N$ with the entries of the vector $\sqrt{N}F_N^{*}[[1+\alpha(-1)^{N}]{\bf 1}_N-\alpha(-1)^Ne_1]$ as its eigenvalues. By a routine calculation, we obtain the diagonalization formula of $\mathcal{H}(B_{\alpha})$ as follows
	\begin{equation}\label{hbepsdiagform}
		\mathcal{H}(B_{\alpha})=\tilde{F}_{N}\Lambda_{\alpha} \tilde{F}_{N}^{*},\quad \Lambda_{\alpha}=N[1+\alpha(-1)^N]{\rm diag}(e_1)-\alpha(-1)^N I_N.
	\end{equation}
	
	By setting $\alpha=0$ in \eqref{hbepsdiagform}, we know that $\mathcal{H}(B)$ is also diagonalizable by $\tilde{F}_N$ as follows
	\begin{equation}\label{hbexpan}
		\mathcal{H}(B)=\tilde{F}_{N}\tilde{\Lambda}\tilde{F}_{N}^{*},\quad \tilde{\Lambda}=N{\rm diag}(e_1).
	\end{equation}

	Let $\sigma(\cdot)$ denotes the spectrum of a square matrix.
	\begin{lemma}\label{pinvkspectrathm}
		The eigenvalues of $\sigma({\bf P}^{-1}{\bf K})$ are lower and upper bounded by $\frac{1}{2}$ and $1$, respectively, i.e., 
		$\sigma({\bf P}^{-1}{\bf K})\subset\left[\frac{1}{2},1\right]$.
	\end{lemma}
	\begin{proof}
		${\bf P}\succ O$ and ${\bf K}\succ O$ imply that  $\sigma({\bf P}^{-1}{\bf K})=\sigma({\bf P}^{-\frac{1}{2}}{\bf K}{\bf P}^{-\frac{1}{2}})\subset (0,+\infty)$ and that eigenvectors of ${\bf P}_{\alpha}^{-1}{\bf P}$ are all real. 
		Let $(\lambda,{\bf x})$ be an eigenpair of ${\bf P}^{-1}{\bf K}$.  Then,
		\begin{equation*}
			{\bf K}{\bf x}=\lambda {\bf P}{\bf x}.
		\end{equation*}
		Hence,
		\begin{equation}\label{pinvkrayleighquotient}
			\lambda=\frac{{\bf x}^{\rm T}{\bf K}{\bf x}}{{\bf x}^{\rm T}{\bf P}{\bf x}}.
		\end{equation}
		The rest of the proof is devoted to estimating the lower and upper bounds of Rayleigh quotient \eqref{pinvkrayleighquotient}.
		From the definitions of ${\bf P}$ and ${\bf K}$, we have
		\begin{align*}
			{\bf K}&=(\sqrt{\tau}I)(\sqrt{\tau}I)+(\sqrt{\eta}G)(\sqrt{\eta}G)^{\rm T},\\
			{\bf P}&=(\sqrt{\tau}I+\sqrt{\eta}G)(\sqrt{\tau}I+\sqrt{\eta}G)^{\rm T}\\
			&=(\sqrt{\eta}G)(\sqrt{\eta}G)^{\rm T}+(\sqrt{\tau}I)(\sqrt{\tau}I)+(\sqrt{\tau}I)(\sqrt{\eta}G)^{\rm T}+(\sqrt{\eta}G)(\sqrt{\tau}I),
		\end{align*}
		where $I=I_N\otimes I_J$ denotes the $(NJ)\times (NJ)$ identity matrix.
		Denote 
		\begin{equation*}
			{\bf c}=\sqrt{\eta}G^{\rm T}{\bf x},\quad {\bf d}=\sqrt{\tau}{\bf x}.
		\end{equation*}
		Then, it is straightforward to verify that
		\begin{equation}\label{rlquotexp}
			\frac{{\bf x}^{\rm T}{\bf K}{\bf x}}{{\bf x}^{\rm T}{\bf P}{\bf x}}=\frac{{\bf c}^{\rm T}{\bf c}+{\bf d}^{\rm T}{\bf d}}{{\bf c}^{\rm T}{\bf c}+{\bf d}^{\rm T}{\bf d}+{\bf c}^{\rm T}{\bf d}+{\bf d}^{\rm T}{\bf c}}
		\end{equation}
		Note that
		\begin{align*}
			{\bf c}^{\rm T}{\bf d}+{\bf d}^{\rm T}{\bf c}=\sqrt{\tau\eta}{\bf x}^{\rm T}(G+ G^{\rm T}){\bf x}=\sqrt{\tau\eta}{\bf x}^{\rm T}[4\mathcal{H}(B)\otimes I_J+2\tau I_N\otimes L_{h}]{\bf x}.
		\end{align*}
		From \eqref{hbexpan}, we know that $\mathcal{H}(B)\succeq O$, which together with $ L_{h}\succeq O$ implies that ${\bf c}^{\rm T}{\bf d}+{\bf d}^{\rm T}{\bf c}\geq 0$. Therefore,
		\begin{equation}\label{eigupbd}
			\frac{{\bf c}^{\rm T}{\bf c}+{\bf d}^{\rm T}{\bf d}}{{\bf c}^{\rm T}{\bf c}+{\bf d}^{\rm T}{\bf d}+{\bf c}^{\rm T}{\bf d}+{\bf d}^{\rm T}{\bf c}}\leq 1.
		\end{equation}
		
		Let $\langle\cdot,\cdot\rangle$ denote the standard inner product in the Euclidean space. Then, by Cauchy-Schwartz inequality, it follows that
		\begin{equation*}
			{\bf c}^{\rm T}{\bf d}+{\bf d}^{\rm T}{\bf c}=2\langle {\bf c},{\bf d}\rangle\leq 2\langle {\bf c},{\bf c}\rangle^{\frac{1}{2}}\langle {\bf d},{\bf d}\rangle^{\frac{1}{2}}\leq \langle {\bf c},{\bf c}\rangle+\langle {\bf d},{\bf d}\rangle={\bf c}^{\rm T}{\bf c}+{\bf d}^{\rm T}{\bf d},
		\end{equation*}
		and thus ${\bf c}^{\rm T}{\bf c}+{\bf d}^{\rm T}{\bf d}+{\bf c}^{\rm T}{\bf d}+{\bf d}^{\rm T}{\bf c}\leq 2({\bf c}^{\rm T}{\bf c}+{\bf d}^{\rm T}{\bf d})$. That means
		\begin{equation*}
			\frac{{\bf c}^{\rm T}{\bf c}+{\bf d}^{\rm T}{\bf d}}{{\bf c}^{\rm T}{\bf c}+{\bf d}^{\rm T}{\bf d}+{\bf c}^{\rm T}{\bf d}+{\bf d}^{\rm T}{\bf c}}\geq \frac{{\bf c}^{\rm T}{\bf c}+{\bf d}^{\rm T}{\bf d}}{2({\bf c}^{\rm T}{\bf c}+{\bf d}^{\rm T}{\bf d})}=\frac{1}{2},
		\end{equation*}
		which together with \eqref{pinvkrayleighquotient}, \eqref{rlquotexp} and \eqref{eigupbd} implies that
		\begin{equation*}
			\frac{1}{2}\leq \lambda\leq 1.
		\end{equation*}
		The proof is complete.
	\end{proof}

	\begin{theorem}\label{pepsinvertibilitythm}
		For any $\alpha\in(0,1]
		\cap\left(0,\frac{\tau}{2\sqrt{\gamma}}\right)$, ${\bf P}_{\alpha}$ is invertible.
	\end{theorem}
	\begin{proof}
		Recall that ${\bf P}_{\alpha}={\bf R}_{\alpha}{\bf R}_{\alpha}^{\rm T}$. Hence, to prove the invertibility of ${\bf P}_{\alpha}$, it suffices to prove the invertibility of ${\bf R}_{\alpha}$. Let $(\lambda,{\bf x})$ be an eigenpair of ${\bf R}_{\alpha}$. Then, ${\bf R}_{\alpha}{\bf x}=\lambda{\bf x}$ and thus
		\begin{align*}
			{\bf x}^{*}{\bf R}_{\alpha}{\bf x}=\lambda||{\bf x}||_2^2\Longrightarrow {\bf x}^{*}\mathcal{H}({\bf R}_{\alpha}){\bf x}=\Re(\lambda)||{\bf x}||_2^2,
		\end{align*}
		where $\Re(\cdot)$ denotes real part of a complex number. If $\mathcal{H}({\bf R}_{\alpha})\succ O$, then $$\Re(\lambda)=\frac{{\bf x}^{*}\mathcal{H}({\bf R}_{\alpha}){\bf x}}{||{\bf x}||_2^2}>0,$$
		which implies $\lambda\neq 0$. The rest  part of this proof is devoted to showing  $\mathcal{H}({\bf R}_{\alpha})\succ O$. 
		
		{\bf Case (i):} $N$ is an odd number. By \eqref{hbepsdiagform} and $\alpha\in(0,1]$, it is clear that $\mathcal{H}(B_{\alpha})\succeq O$. Combining this with \eqref{hrepsexpan}, we know that $\mathcal{H}({\bf R}_{\alpha})\succ O$ in this case.
		
		{\bf Case (ii):} $N$ is an even number. Denote $\tilde{\Lambda}_{\alpha}=\sqrt{\tau}I_N+2\sqrt{\eta}\Lambda_{\alpha}$. From \eqref{hbepsdiagform} and \eqref{hrepsexpan}, we know that
		\begin{equation}\label{hrepsprelb}
			\mathcal{H}({\bf R}_{\alpha})\succeq (\sqrt{\tau}I_N+2\sqrt{\eta}\mathcal{H}(B_{\alpha}))\otimes I_J=(\tilde{F}_{N}\otimes I_J)(\tilde{\Lambda}_{\alpha}\otimes I_J)(\tilde{F}_{N}^{*}\otimes I_J).
		\end{equation}
		Hence, to show $\mathcal{H}({\bf R}_{\alpha})\succ O$, it suffices to show the diagonal entries of the diagonal matrix $\tilde{\Lambda}_{\alpha}$ are all positive. Since $\alpha\in(0,1]$, $\tilde{\Lambda}_{\alpha}(1,1)=\sqrt{\tau}+N[1+\alpha(-1)^N]-\alpha(-1)^N=\sqrt{\tau}+N+(N-1)\alpha>0.$  Moreover, since $\alpha\in\left(0,\frac{\tau}{2\sqrt{\gamma}}\right)$, $$\tilde{\Lambda}_{\alpha}(i,i)=\sqrt{\tau}-2\sqrt{\eta}\alpha> \sqrt{\tau}-2\sqrt{\eta}\tau/(2\sqrt{\gamma})=0.$$ Hence, $\mathcal{H}({\bf R}_{\alpha})\succ O$.
		
		To summarize the discussion above, we see that $\lambda\neq 0$. By the generality of $\lambda$, we know that ${\bf R}_{\alpha}$ is invertible and thus ${\bf P}_{\alpha}$ is invertible. The proof is complete.
	\end{proof}
	
	Let $\rho(\cdot)$ denote the spectral radius of a square matrix.
	\begin{lemma}\label{resmatsesti}
		\begin{align*}
			||\tilde{B}\tilde{B}^{\rm T}||_2\leq \frac{2T^2}{\tau^2}.
		\end{align*}
	\end{lemma}
	\begin{proof}
		A routine calculation yields that
		\begin{equation*}
			(\tilde{B}\tilde{B}^{\rm T})(i,j)=\begin{cases}
				4(-1)^{j-i}(N-j),\quad j\geq i,\\
				(\tilde{B}\tilde{B}^{\rm T})(j,i),\quad j<i.
			\end{cases}
		\end{equation*}
		Because of the equality, it is easy to check that
		\begin{align*}
			||\tilde{B}\tilde{B}^{\rm T}||_1=||[\tilde{B}\tilde{B}^{\rm T}](:,1)||_1=4\sum\limits_{i=1}^{N}(N-i)=2N(N-1)<2N^2=\frac{2T^2}{\tau^2},
		\end{align*} 
		where $[\tilde{B}\tilde{B}^{\rm T}](:,1)$ denotes the first column of $\tilde{B}\tilde{B}^{\rm T}$.
		Then, by symmetry of $\tilde{B}\tilde{B}^{\rm T}$, we know that
		\begin{equation*}
			||\tilde{B}\tilde{B}^{\rm T}||_2=\rho(\tilde{B}\tilde{B}^{\rm T})\leq ||\tilde{B}\tilde{B}^{\rm T}||_1<\frac{2T^2}{\tau^2}.
		\end{equation*} 
		The proof is complete.
	\end{proof}

	\begin{lemma}\label{invpepspspectralm}
		For any $$\alpha\in\Bigg(0,\min\left\{\frac{\tau}{24\sqrt{\gamma}},\frac{\tau^{\frac{3}{2}}}{2\sqrt{6\gamma}T},\frac{\tau^{2}}{8\sqrt{3\gamma}T},\frac{1}{3}\right\}\Bigg],$$ it holds that $\sigma({\bf P}_{\alpha}^{-1}{\bf P})\subset \left[\frac{3}{4},\frac{3}{2}\right]$.
	\end{lemma}
	\begin{proof}
		Since $\alpha\in\Big(0,\frac{1}{3}\Big]\cap\Big(0,\frac{\tau}{24\sqrt{\gamma}}\Big]\subset(0,1]\cap\left(0,\frac{\tau}{2\sqrt{\gamma}}\right)$, Lemma \ref{pepsinvertibilitythm}${\bf (i)}$ implies that ${\bf P}_{\alpha}^{-1}$ exists. Moreover,  ${\bf P}_{\alpha}\succ O$ and ${\bf P}\succ O$ imply that  $\sigma({\bf P}_{\alpha}^{-1}{\bf P})=\sigma({\bf P}_{\alpha}^{-\frac{1}{2}}{\bf P}{\bf P}_{\alpha}^{-\frac{1}{2}})\subset (0,+\infty)$ and that eigenvectors of ${\bf P}_{\alpha}^{-1}{\bf P}$ are all real. Let $(\lambda,{\bf x})$ be an eigenpair of ${\bf P}_{\alpha}^{-1}{\bf P}$. Then, ${\bf P}{\bf x}=\lambda{\bf P}_{\alpha}{\bf x}$ and thus 
		\begin{equation}\label{pepsinvpeigexpress}
			\lambda=\frac{{\bf x}^{\rm T}{\bf P}{\bf x}}{{\bf x}^{\rm T}{\bf P}_{\alpha}{\bf x}}.
		\end{equation}
		From definition of ${\bf P}_{\alpha}$ and ${\bf P}$, we know that they can be written as 
		\begin{align}
			{\bf P}=&[\tau I_N+4\sqrt{\gamma}\mathcal{H}(B)+4\eta BB^{\rm T}]\otimes I_J+[2\tau^{\frac{3}{2}}\eta^{\frac{1}{2}}I_N+4\gamma  \mathcal{H}(B)]\otimes L_{h}+\tau^2\eta I_N\otimes L_{h}^2,\label{pexpan}\\
			{\bf P}_{\alpha}=&[\tau I_N+4\sqrt{\gamma}\mathcal{H}(B_{\alpha})+4\eta BB^{\rm T}]\otimes I_J+[2\tau^{\frac{3}{2}}\eta^{\frac{1}{2}}I_N+4\gamma \mathcal{H}(B_{\alpha})]\otimes L_{h}+\tau^2\eta I_N\otimes L_{h}^2\notag\\
			&+4\eta\alpha(B\tilde{B}^{\rm T}+\tilde{B}B^{\rm T})\otimes I_J+4\eta\alpha^2(\tilde{B}\tilde{B}^{\rm T})\otimes I_J.\label{pepsexpan}
		\end{align}
		The following inequalities
		\begin{align*}
			&O\preceq \left(\frac{1}{\sqrt{3}} B+\sqrt{3}\alpha\tilde{B}\right)\left(\frac{1}{\sqrt{3}} B+\sqrt{3}\alpha\tilde{B}\right)^{\rm T}=\alpha(\tilde{B}B^{\rm T}+B\tilde{B}^{\rm T})+\frac{1}{3} BB^{\rm T}+3\alpha^2\tilde{B}\tilde{B}^{\rm T},\\
			&O\preceq \left(\frac{1}{\sqrt{3}} B-\sqrt{3}\alpha\tilde{B}\right)\left(\frac{1}{\sqrt{3}} B-\sqrt{3}\alpha\tilde{B}\right)^{\rm T}=-\alpha(\tilde{B}B^{\rm T}+B\tilde{B}^{\rm T})+\frac{1}{3} BB^{\rm T}+3\alpha^2\tilde{B}\tilde{B}^{\rm T},
		\end{align*}
		imply that
		\begin{equation}\label{medresesti}
			-\frac{1}{3} BB^{\rm T}-3\alpha^2\tilde{B}\tilde{B}^{\rm T}\preceq	\alpha(B\tilde{B}^{\rm T}+\tilde{B}B^{\rm T})	\preceq \frac{1}{3} BB^{\rm T}+3\alpha^2\tilde{B}\tilde{B}^{\rm T}.
		\end{equation}
		
		Moreover, Lemma \ref{resmatsesti} and $0<\alpha\leq \frac{\tau^{2}}{8\sqrt{3\gamma}T}\leq \frac{\tau^2}{4\sqrt{6\gamma}T}$ imply that
		\begin{align}
			-8\eta\alpha^2\tilde{B}\tilde{B}^{\rm T}&\succeq \frac{-16\eta T^2\alpha^2}{\tau^2}I_N\notag\\
			&\succeq \left[ \frac{-16\eta T^2}{\tau^2}\times \left(\frac{\tau^2}{4\sqrt{6\gamma}T}\right)^2\right]I_N=\frac{-\tau}{6} I_N.\label{tildbsqlbesti1}
		\end{align}
		Then,
		\begin{align}
			&\tau I_N+4\sqrt{\gamma}\mathcal{H}(B_{\alpha})+4\eta BB^{\rm T}+4\eta\alpha(B\tilde{B}^{\rm T}+\tilde{B}B^{\rm T})+4\eta\alpha^2\tilde{B}\tilde{B}^{\rm T}\notag\\
			&\succeq \tau I_N+4\sqrt{\gamma}\mathcal{H}(B_{\alpha})+\frac{8}{3}\eta BB^{\rm T}-8\eta\alpha^2\tilde{B}\tilde{B}^{\rm T}\notag\\
			&=(\tau-4\sqrt{\gamma}(-1)^{N}\alpha)I_N+4\sqrt{\gamma}[1+\alpha(-1)^N]\mathcal{H}(B)+\frac{8}{3}\eta BB^{\rm T}-8\eta\alpha^2\tilde{B}\tilde{B}^{\rm T}\notag\\
			&\succeq \left[\tau-4\sqrt{\gamma}(-1)^{N}\alpha-\tau/6\right]I_N+4\sqrt{\gamma}[1+\alpha(-1)^N]\mathcal{H}(B)+\frac{8}{3}\eta BB^{\rm T}\notag\\
			&\succeq \frac{2}{3}\tau I_N+\frac{8}{3}\sqrt{\gamma}\mathcal{H}(B)+\frac{8}{3}\eta BB^{\rm T},\label{pepspart1lbesti}
		\end{align}
		where the first inequality comes from first inequality of \eqref{medresesti}, the first equality comes from \eqref{hbepsdiagform} and \eqref{hbexpan}, the second inequality comes from \eqref{tildbsqlbesti1}, the third inequality comes from $\alpha\in\Big(0,\frac{1}{3}\Big]\cap\Big(0,\frac{\tau}{24\sqrt{\gamma}}\Big]$.
		
		By \eqref{hbepsdiagform}, \eqref{hbexpan} and $\alpha\in \Big(0,\frac{\tau}{24\sqrt{\gamma}}\Big]\cap\Big(0,\frac{1}{3}\Big]$, we obtain
		\begin{align*}
			2\tau^{\frac{3}{2}}\eta^{\frac{1}{2}}I_N+4\gamma \mathcal{H}(B_{\alpha})&=[2\tau^{\frac{3}{2}}\eta^{\frac{1}{2}}-4\gamma\alpha(-1)^N ]I_N+4\gamma[1+\alpha(-1)^N]\mathcal{H}(B)\notag\\
			&\succeq \left[2\tau^{\frac{3}{2}}\eta^{\frac{1}{2}}-4\gamma\times\frac{\tau}{24\sqrt{\gamma}}\right]I_N +\frac{8}{3}\gamma\mathcal{H}(B)\notag\\
			&=\frac{11}{6}\tau^{\frac{3}{2}}\eta^{\frac{1}{2}}I_N+\frac{8}{3}\gamma\mathcal{H}(B),
		\end{align*}
		which together with \eqref{pepsexpan}, \eqref{pepspart1lbesti} and \eqref{pexpan} implies that
		\begin{equation*}
			{\bf P}_{\alpha}\succeq \frac{2}{3}{\bf P}.
		\end{equation*}
		Applying the above inequality to \eqref{pepsinvpeigexpress}, we obtain that  
		\begin{equation}\label{eigubodd}
			\lambda\leq\frac{3}{2}.
		\end{equation}
		
		By Lemma \ref{resmatsesti} and $0<\alpha\leq \frac{\tau^2}{8\sqrt{3\gamma}T}$, we obtain
		\begin{align}
			16\eta\alpha^2\tilde{B}\tilde{B}^{\rm T}&\preceq \frac{32\eta T^2\alpha^2}{\tau^2}I_N\notag\\
			&\preceq  \left(\frac{\tau^2}{8\sqrt{3\gamma}T}\right)^2\times \frac{32\eta T^2}{\tau^2}I_N =\frac{\tau}{6}I_N,\notag
		\end{align}
		which together with \eqref{medresesti}, \eqref{hbepsdiagform}, \eqref{hbexpan} and $\alpha\in\Big(0,\frac{1}{3}\Big]\cap\Big(0,\frac{\tau}{24\sqrt{\gamma}}\Big]$ implies that
		\begin{align}
			&\tau I_N+4\sqrt{\gamma}\mathcal{H}(B_{\alpha})+4\eta BB^{\rm T}+4\eta\alpha(B\tilde{B}^{\rm T}+\tilde{B}B^{\rm T})+4\eta\alpha^2\tilde{B}\tilde{B}^{\rm T}\notag\\
			&\preceq \tau I_N+4\sqrt{\gamma}\mathcal{H}(B_{\alpha})+\frac{16}{3}\eta BB^{\rm T}+16\eta\alpha^2\tilde{B}\tilde{B}^{\rm T}\notag\\
			&=(\tau-4\sqrt{\gamma}(-1)^{N}\alpha)I_N+4\sqrt{\gamma}[1+\alpha(-1)^N]\mathcal{H}(B)+ \frac{16\eta}{3}BB^{\rm T}+16\eta\alpha^2\tilde{B}\tilde{B}^{\rm T}\notag\\
			&\preceq \left[\tau-4\sqrt{\gamma}(-1)^{N}\alpha+\tau/6\right]I_N+4\sqrt{\gamma}[1+\alpha(-1)^N]\mathcal{H}(B)+\frac{16}{3}\eta BB^{\rm T}\notag\\
			&\preceq \frac{4\tau}{3}I_N+\frac{16\gamma}{3}\sqrt{\gamma}\mathcal{H}(B)+\frac{16\eta}{3}BB^{\rm T}.\label{pepspart1ubesti}
		\end{align}
		
		By \eqref{hbepsdiagform}, \eqref{hbexpan} and $\alpha\in \Big(0,\frac{\tau}{24\sqrt{\gamma}}\Big]\cap(0,\frac{1}{3}]$, we obtain
		\begin{align*}
			2\tau^{\frac{3}{2}}\eta^{\frac{1}{2}}I_N+4\gamma \mathcal{H}(B_{\alpha})&=[2\tau^{\frac{3}{2}}\eta^{\frac{1}{2}}-4\gamma\alpha(-1)^N ]I_N+4\gamma[1+\alpha(-1)^N]\mathcal{H}(B)\notag\\
			&\preceq \left[2\tau^{\frac{3}{2}}\eta^{\frac{1}{2}}+4\gamma\times\frac{\tau}{24\sqrt{\gamma}}\right]I_N +\frac{16\gamma}{3}\mathcal{H}(B)\notag\\
			&=\frac{13}{6}\tau^{\frac{3}{2}}\eta^{\frac{1}{2}}I_N+\frac{16\gamma}{3}\mathcal{H}(B),
		\end{align*}
		which together with \eqref{pepsexpan}, \eqref{pepspart1ubesti} and \eqref{pexpan} implies that
		\begin{equation*}
			{\bf P}_{\alpha}\preceq \frac{4}{3}{\bf P}.
		\end{equation*}
		Applying the above inequality to \eqref{pepsinvpeigexpress}, we obtain that
		\begin{equation*}
			\lambda\geq \frac{3}{4},
		\end{equation*}
		which together with \eqref{eigubodd} completes the proof.
	\end{proof}
	
	With the help of Lemma \ref{pinvkspectrathm} and Lemma \ref{invpepspspectralm}, we obtain an estimation of range of the spectrum of the preconditioned matrix ${\bf P}_{\alpha}^{-1}{\bf K}$ in the following theorem.
	\begin{theorem}\label{precedmatspectrathm}
		For any $$\alpha\in\Bigg(0,\min\left\{\frac{\tau}{24\sqrt{\gamma}},\frac{\tau^{\frac{3}{2}}}{2\sqrt{6\gamma}T},\frac{\tau^{2}}{8\sqrt{3\gamma}T},\frac{1}{3}\right\}\Bigg],$$  it holds that $\sigma({\bf P}_{\alpha}^{-1}{\bf K})\subset \left[\frac{3}{8},\frac{3}{2}\right]$.
	\end{theorem}
	\begin{proof}
		Since ${\bf P}_{\alpha}^{-1}{\bf K}$ is similar to ${\bf P}_{\alpha}^{-\frac{1}{2}}{\bf K}{\bf P}_{\alpha}^{-\frac{1}{2}}$, we see that ${\bf P}_{\alpha}^{-1}{\bf K}$ has real eigenvalues and real eigenvectors. Let $(\lambda,{\bf x})$ be an eigenpair of ${\bf P}_{\alpha}^{-1}{\bf K}$. Then, ${\bf K}{\bf x}=\lambda{\bf P}_{\alpha}{\bf x}$ and thus
		\begin{equation*}
			\lambda=\frac{{\bf x}^{\rm T}{\bf K}{\bf x}}{{\bf x}^{\rm T}{\bf P}_{\alpha}{\bf x}}=\frac{{\bf x}^{\rm T}{\bf K}{\bf x}}{{\bf x}^{\rm T}{\bf P}{\bf x}}\times \frac{{\bf x}^{\rm T}{\bf P}{\bf x}}{{\bf x}^{\rm T}{\bf P}_{\alpha}{\bf x}}.
		\end{equation*}
		Then, Lemma \ref{pinvkspectrathm} and Lemma \ref{invpepspspectralm} imply that
		\begin{align*}
			\frac{3}{8}\leq \lambda\leq \frac{3}{2},
		\end{align*}
		which completes the proof.
	\end{proof}
	
	\begin{remark}\label{looseboundremark}
		Actually, by the fact that $\lim\limits_{\alpha\rightarrow 0^{+}}{\bf P}_{\alpha}={\bf P}$, one can show that for sufficiently small $\alpha$, ${\bf P}_{\alpha}^{-1}{\bf P}$ is close to the identity matrix and thus the interval $[0.5,1]$ covering $\sigma({\bf P}^{-1}{\bf K})$ almost covers  $\sigma({\bf P}_{\alpha}^{-1}{\bf K})$. In other words, one can show that $\sigma({\bf P}_{\alpha}^{-1}{\bf K})\subset \left[0.5-\epsilon(\alpha),1+\epsilon(\alpha)\right]$ for some nonnegative function $\epsilon(\cdot)$ with $\lim\limits_{\alpha\rightarrow 0^{+}}\epsilon(\alpha)=0$ as a sharper result to Theorem \ref{precedmatspectrathm}. The reason why we give a loose interval $\left[\frac{3}{8},\frac{3}{2}\right]$  covering $\sigma({\bf P}_{\alpha}^{-1}{\bf K})$ is to simplify the expression of the upper bound for $\alpha$ and to give a tidy convergence factor of the PCG solver in Theorem \ref{pcgcvgthm}.
	\end{remark}
	
	For any symmetric positive definite matrix $A\in\mathbb{R}^{n\times n}$, one can define an inner product
	\begin{equation*}
		(x,z)_{A}:=x^{\rm T}Az,\quad x,z\in\mathbb{R}^{n\times 1},
	\end{equation*}
	which induces a vector norm as follows
	\begin{equation*}
		||x||_{A}=(x,x)_{A}^{\frac{1}{2}}.
	\end{equation*}
	As ${\bf K}\succ O$, one can define inner product as follows
	\begin{equation*}
		(x,z)_{{\bf K}}:=x^{\rm T}{\bf K}z,\quad x,z\in\mathbb{R}^{NJ\times 1}.
	\end{equation*}
	Correspondingly, one can define the norm induced by the inner product as follows
	\begin{equation*}
		||x||_{{\bf K}}:=(x,x)_{{\bf K}}^{\frac{1}{2}},\quad x\in\mathbb{R}^{n\times 1}.
	\end{equation*}
	\begin{lemma}\textnormal{(see \cite{axelsson1986rate})}\label{pcgcvglm}
		Let $A_1,A_2\in\mathbb{R}^{n\times n}$ be symmetric positive definite matrices. Then, the convergence rate of the PCG solver for the preconditioned system $A_1^{-1}A_2{\bf x}=A_1^{-1}{\bf b}$ can be estimated as follows
		\begin{equation*}
			||{\bf x}^{k}-{\bf x}||_{A_2}\leq 2\left(\frac{\sqrt{\hat{c}}-\sqrt{\check{c}}}{\sqrt{\hat{c}}+\sqrt{\check{c}}}\right)^{k}||{\bf x}^{0}-{\bf x}||_{A_2},
		\end{equation*}
		where $\hat{c}$, $\check{c}$ are any positive numbers satisfying $\sigma(A_1^{-1}A_2)\subset[\check{c},\hat{c}]$, ${\bf x}^{k}$ $(k\geq 1)$ denotes the $k$th iterative solution by PCG solver, ${\bf x}^{0}$ denotes an arbitrary initial guess.
	\end{lemma}
	
	Applying Lemma \ref{pcgcvglm} to Theorem \ref{precedmatspectrathm}, we immediately obtain the  convergence result of PCG solver for our preconditioned system \eqref{precschurcmpsys} as stated in the following theorem.
	\begin{theorem}\label{pcgcvgthm}
		Take $$\alpha\in (0,\nu],$$ with
		\begin{equation*}
			\nu=\min\left\{\frac{\tau}{24\sqrt{\gamma}},\frac{\tau^{\frac{3}{2}}}{2\sqrt{6\gamma}T},\frac{\tau^{2}}{8\sqrt{3\gamma}T},\frac{1}{3}\right\}.
		\end{equation*}
		Then, the PCG solver for the preconditioned system \eqref{precschurcmpsys} has a linear convergence rate independent of matrix size $NJ$, i.e.,
		\begin{equation*}
			||v^{k}-v||_{{\bf K}}\leq \frac{2}{3^k}||v^{0}-v||_{{\bf K}},
		\end{equation*}
		where $v^{k}$ $(k\geq 1)$ denotes the $k$th iterative solution by PCG solver, $v^{0}$ denotes an arbitrary initial guess.
	\end{theorem}
	
	\begin{remark}
		Roughly speaking, in Theorem \ref{pcgcvgthm}, the upper bound for $\alpha$ is of $\mathcal{O}\left(\frac{\tau^2}{\sqrt{\gamma}}\right)$. That means when  $\alpha\lesssim \frac{\tau^2}{\sqrt{\gamma}}$, the PCG solver for the preconditioned system \eqref{precschurcmpsys} has a linear convergence rate independent of matrix size and regularization parameter $\gamma$, which implies that the convergence is optimal. 
	\end{remark}

	\section{Numerical Results}
	In this section, we test the performance of the proposed PinT preconditioner by several numerical examples. All numerical experiments presented in this section are conducted via MATLAB R2018b on a PC with the configuration:
	Intel(R) Core(TM) i7-4720HQ CPU 2.60 GHz and 8 GB RAM.

	The stopping criterion for PCG solver is set as $||{\bf r}_k||_2\leq 10^{-8}||{\bf r}_0||_2$, where ${\bf r}_k=b-{\bf K}v_k$ denotes the unpreconditioned residual vector at $k$th PCG iteration and $v_0$ denotes the zero initial guess. According to Theorem \ref{pcgcvgthm}, we take the following default value for $\alpha$ in this section, if not specified.
	\begin{equation*}
		\alpha=\frac{1}{2}\nu,
	\end{equation*}
	where $\nu$ is given in Theorem \ref{pcgcvgthm}.
	It is easy to see that the so evaluated $\alpha$ satisfies the assumption presented in Theorem \ref{pcgcvgthm}. Central difference scheme on uniform square grids is adopted to discretize the spatial operator $ L=\nabla(a(x)\nabla\cdot)$. When $a(x)$ is a constant function, the linear systems appearing in step-(b) of \eqref{threestepimplemt} are diagonalizable by fast sine transform and solving each of those linear system requires $\mathcal{O}(J\log J)$ flops. For more general diffusion coefficient $a(x)$, we adopt one iteration of V-cycle geometric multigrid method to solve the linear systems in step-(b) of \eqref{threestepimplemt} and the complexity of the V-cycle iteration is optimal (i.e., proportional to the number of unknowns). 
	
	As we mentioned in the introduction part, the authors in \cite{pearson2012regularization} also proposed a preconditioner for Schur complement of matrix arising from optimal control problem. Following the idea in \cite{pearson2012regularization}, we find that the preconditioner proposed in \cite{pearson2012regularization}  is exactly the matrix ${\bf P}$ defined in \eqref{Pmeddef}. We will compare the performance of ${\bf P}$ and our proposed preconditioner ${\bf P}_{\alpha}$ in the following experiment. When ${\bf P}$ is used as a preconditioner, it requires to invert ${\bf P}$ during each PCG-${\bf P}$ iteration. Note that ${\bf P}={\bf R}{\bf R}^{\rm T}$ and $\mathbb{R}$ is a block lower triangular matrix. Hence, the inversion of ${\bf P}$ resorts to block forward and backward substitutions, which is not PinT.

	Recall that the purpose of solving the Schur complement system \eqref{origschurcmpsys} (or equivalently \eqref{precschurcmpsys}) is to compute $[p_{\tau,h}^{\rm T},y_{\tau,h}^{\rm T}]^{\rm T}$ in \eqref{kktinversform}. Let $p_{\tau,h}^{*}$ ($y_{\tau,h}^{*}$)  denote the approximate solution to $p_{\tau,h}$ ($y_{\tau,h}$). Then, we define the error measure as
	\begin{equation*}
		{\rm E}_{N,J}=\left|\left|\left[\begin{array}[c]{c}
			p_{\tau,h}^{*}\\
			y_{\tau,h}^{*}
		\end{array}\right]-\left[\begin{array}[c]{c}
			p_{\tau,h}\\
			y_{\tau,h}
		\end{array}\right]\right|\right|_{\infty}.
	\end{equation*}
	Denote by `Iter', the iteration number of the PCG solver. Denote by `CPU', the computational time in unit of second. Denote by PCG-${\bf P}_{\alpha}$ (PCG-${\bf P}$) the PCG solver with preconditioner ${\bf P}_{\alpha}$ (${\bf P}$) for the Schur complement system. 
	
	\begin{example}\label{constcoeffanalyticalsolexpl1d}
		{\rm 
			In this example, we consider the minimization problem \eqref{optcontrlprob} with
			\begin{align*}
				&a(x)\equiv 1,\quad\Omega=(0,1),\quad T=1,\quad f(x,t)=(\pi^2-1)\sin(\pi x)\exp(-t),\\
				&g(x_1,x_2,t)=\sin(\pi x_2)\sin(\pi x_1)\exp(-t),
			\end{align*}
			the analytical solution of which is given by
			\begin{align*}
				y(x_1,x_2,t)=\sin(\pi x_1)\sin(\pi x_2)\exp(-t),\quad u\equiv 0.
			\end{align*}
			To support the results of Theorem \ref{precedmatspectrathm}, we plot the spectrum of the preconditioned matrices corresponding to different values of $\alpha$ in Figure \ref{spectrumplot}.
			\begin{figure}[H]
				\centering
				\subfigure[$\alpha$=5.0e-1$>\nu\approx$7.22e-2]{\includegraphics[width=1.5in]{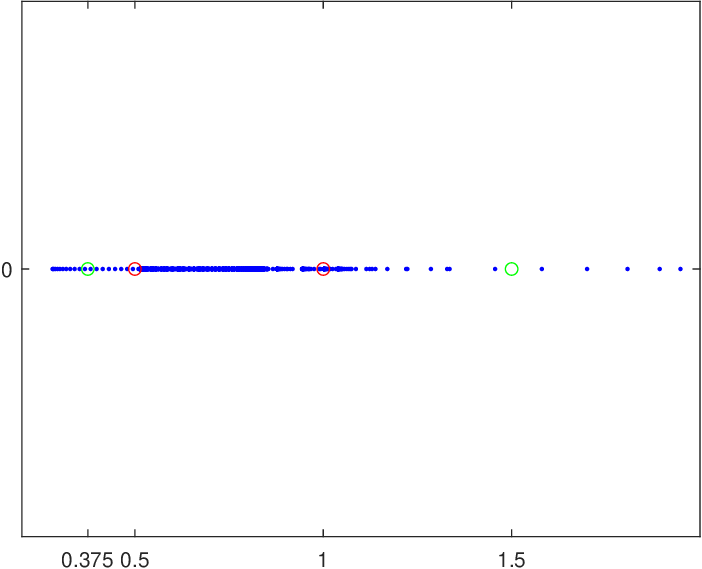}}\quad\subfigure[$\alpha=\nu\approx$7.22e-2]{\includegraphics[width=1.5in]{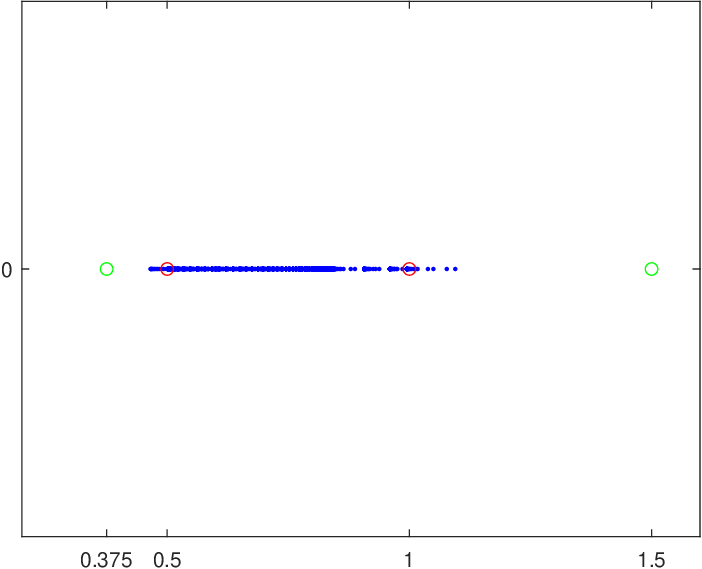}}\quad\subfigure[$\alpha$=5.0e-3$<\nu\approx$7.22e-2]{\includegraphics[width=1.5in]{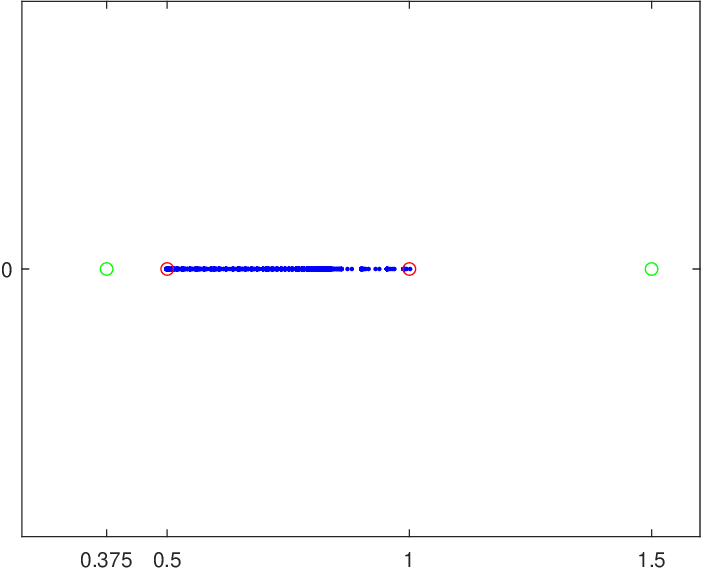}}
				
				\subfigure[$\alpha$=5.0e-4$<\nu\approx$7.22e-2]{\includegraphics[width=1.5in]{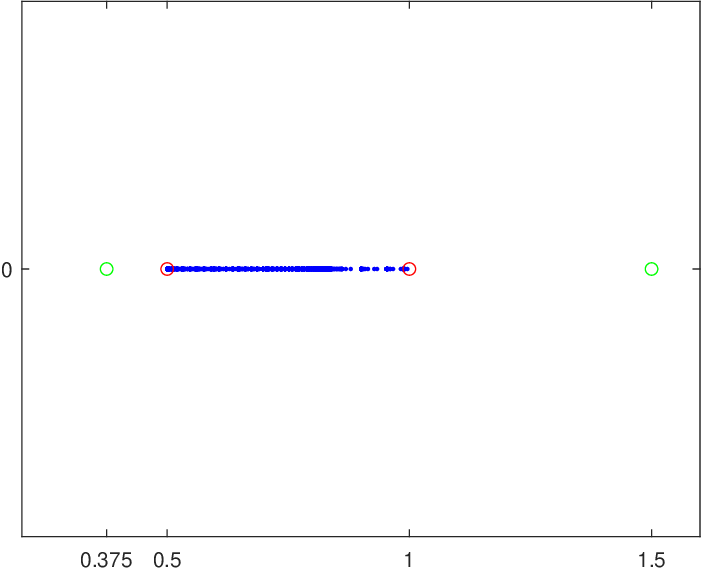}}\quad\subfigure[$\alpha$=5.0e-5$<\nu\approx$7.22e-2]{\includegraphics[width=1.5in]{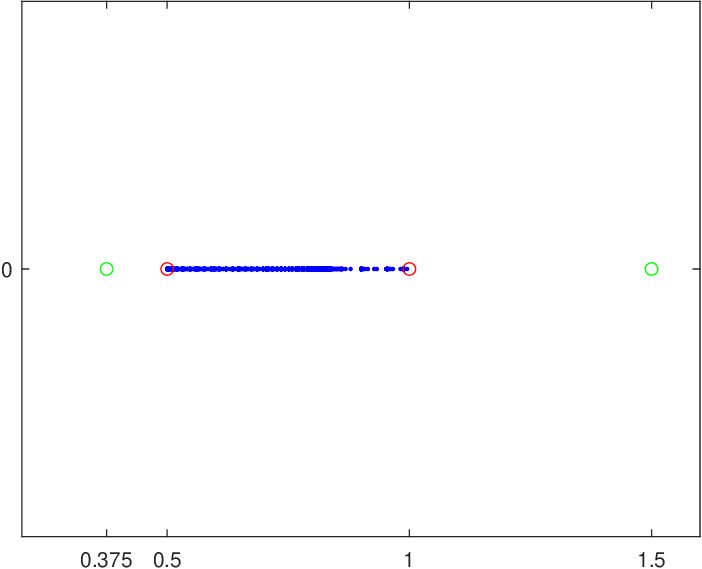}}\quad\subfigure[$\alpha$=5.0e-6$<\nu\approx$7.22e-2]{\includegraphics[width=1.5in]{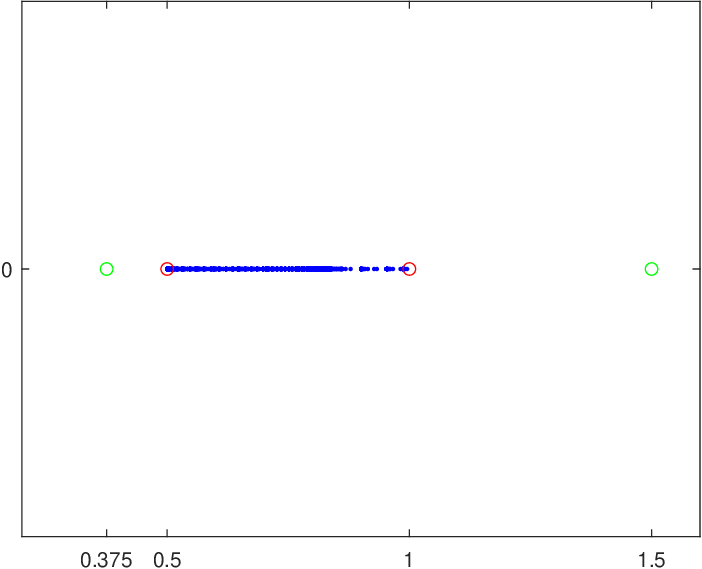}}
				\caption{Spectrum of ${\bf P}_{\alpha}^{-1}{\bf K}$ corresponding to different values of $\alpha$ when $N=40$, $J=63$, $\gamma=\tau^4$. The blue points denotes the eigenvalues of ${\bf P}_{\alpha}^{-1}{\bf K}$. The  red circles denote the two points $x=0.5$ and $x=1$. The green circles denote the two points $x=3/8$ and $x=3/2$.}\label{spectrumplot}
			\end{figure}
			From Figure \ref{spectrumplot}(a), we see that the spectrum of ${\bf P}_{\alpha}^{-1}{\bf K}$ may lie outside the interval $[3/8,3/2]$ when $\alpha>\nu$. That means the setting $\alpha\in(0,\nu]$ proposed in Theorem \ref{pcgcvgthm} makes sense. Figures \ref{spectrumplot}(a)-(b) illustrate that introducing the block $\alpha$-circulant approximations to the classical MSC preconditioner leads to perturbations in eigenvalues so that the spectrum of the preconditioned matrix is no longer completely contained in $[1/2,1]$ as the classical MSC preconditioner behaves. Figures \ref{spectrumplot} (b)-(f) demonstrates that $\sigma({\bf P}_{\alpha}^{-1}{\bf K})$ indeed lies in [3/8,3/2] for $\alpha\in(0,\nu]$, which supports the result of Theorem \ref{pcgcvgthm}. Figures \ref{spectrumplot} (c)-(f) shows that $\sigma({\bf P}_{\alpha}^{-1}{\bf K})$ is almost contained in [1/2,1] as $\sigma({\bf P}^{-1}{\bf K})$ distributes. This is not surprising, since $\lim\limits_{\alpha\rightarrow 0^{+}}{\bf P}_{\alpha}^{-1}{\bf K}={\bf P}^{-1}{\bf K}$. On the other hand, Figures \ref{spectrumplot} (c)-(f) also indicate that [3/8,3/2] is a loose interval covering $\sigma({\bf P}_{\alpha}^{-1}{\bf K})$. The reason why we choose a loose interval has been explained in Remark \ref{looseboundremark}, which is only for a simple bound $\nu$ for $\alpha$ and a tidy convergence rate of the PCG solver  shown in Theorem \ref{pcgcvgthm}.
		}	
	\end{example}

	\begin{example}\label{constcoeffanalyticalsolexpl}
		{\rm 
			In this example, we consider the minimization problem \eqref{optcontrlprob} with
			\begin{align*}
				&a(x_1,x_2)\equiv 1,\quad\Omega=(0,1)^2,\quad T=1,\quad f(x_1,x_2,t)=(2\pi^2-1)\sin(\pi x_2)\sin(\pi x_1)\exp(-t),\\
				&g(x_1,x_2,t)=\sin(\pi x_2)\sin(\pi x_1)\exp(-t),
			\end{align*}
			the analytical solution of which is given by
			\begin{align*}
				y(x_1,x_2,t)=\sin(\pi x_1)\sin(\pi x_2)\exp(-t),\quad u\equiv 0.
			\end{align*}
			The numerical results of PCG-${\bf P}_{\alpha}$ for solving Example \ref{constcoeffanalyticalsolexpl} are presented in Table \ref{explconstanalyticaltab}. We observe from Table \ref{explconstanalyticaltab} that (i) the iteration numbers of the two solvers in Table \ref{explconstanalyticaltab} changes slightly as $N$ or $J$ changes; (ii) the iteration numbers of the two solvers are bounded (increasing first and then decreasing) with respect to the changes of $\gamma$. The boundedness of iteration number of PCG-${\bf P}_{\alpha}$ in Table \ref{explconstanalyticaltab} with respect to $N$, $J$ and $\gamma$ supports the theoretical results in Theorem \ref{pcgcvgthm}. Another interesting observation from Table \ref{constcoeffanalyticalsolexpl} is that the iteration numbers, errors of the two solvers are roughly the same while the CPU cost of PCG-${\bf P}_{\alpha}$ is much smaller than that of PCG-${\bf P}$ especially for large $N$. This is because ${\bf P}_{\alpha}$ is a PinT preconditioner while ${\bf P}$ is not. The smaller CPU cost of PCG-${\bf P}_{\alpha}$ compared with PCG-${\bf P}$ demonstrates the advantage of a PinT preconditioner.  
			\newpage
			\begin{table}[H]
				\begin{center}
					\caption{Performance of the PCG-${\bf P}_{\alpha}$ and PCG-${\bf P}$ on Example \ref{constcoeffanalyticalsolexpl}.
					}\label{explconstanalyticaltab}
					\setlength{\tabcolsep}{0.7em}
					\begin{tabular}[c]{cccc|ccc|ccc}
						\hline
						\multirow{2}{*}{$\gamma$} &\multirow{2}{*}{$N$} &\multirow{2}{*}{$J$}&\multirow{2}{*}{$\alpha$}&\multicolumn{3}{c|}{PCG-${\bf P}_{\alpha}$}&\multicolumn{3}{c}{PCG-${\bf P}$}\\
						\cline{5-10}
						&&&&Iter & CPU(s)& ${\rm E}_{N,J}$&Iter & CPU(s)& ${\rm E}_{N,J}$\\
						\hline
						\multirow{9}{*}{1e-7}&\multirow{3}{*}{200}&961   &\multirow{3}{*}{2.85e-3}   &4 &0.66 &4.43e-3 &4 &4.04 &4.43e-3\\
						&                    &3969  &                           &4 &3.06 &4.43e-3 &4 &12.36&4.43e-3\\
						&                    &16129 &                           &4 &11.62&4.43e-3 &4 &42.77&4.44e-3\\
						\cline{2-10}
						&\multirow{3}{*}{400}&961   &\multirow{3}{*}{7.13e-4}   &4 &1.36 &1.99e-3 &4 &12.20 &1.99e-3\\
						&                    &3969  &                           &4 &5.66 &1.99e-3 &4 &41.32 &1.99e-3\\
						&                    &16129 &                           &4 &23.38&1.99e-3 &4 &137.84&1.99e-3\\
						\cline{2-10}
						&\multirow{3}{*}{800}&961   &\multirow{3}{*}{1.78e-4}   &4 &2.79 &8.29e-4 &4 &37.75 &8.29e-4\\
						&                    &3969  &                           &4 &11.35&8.29e-4 &4 &142.26&8.29e-4\\
						&                    &16129 &                           &4 &46.81&8.29e-4 &4 &446.34&8.29e-4\\		
						\hline
						\multirow{9}{*}{1e-5}&\multirow{3}{*}{200}&961   &\multirow{3}{*}{2.85e-4}   &6 &0.90 &2.45e-3 &6 &5.17 &2.4e-3\\
						&                    &3969  &   &6 &3.77 &2.45e-3 &6 &15.46&2.45e-3\\
						&                    &16129 &  &6 &15.70&2.45e-3 &6 &53.49&2.45e-3\\
						\cline{2-10}
						&\multirow{3}{*}{400}&961   &\multirow{3}{*}{7.13e-5}    &7 &1.98 &1.22e-3 &6 &15.53 &1.22e-3\\
						&                    &3969  &                            &7 &8.81 &1.22e-3 &6 &50.17 &1.22e-3\\
						&                    &16129 &                            &7 &35.41&1.22e-3 &6 &179.41&1.22e-3\\
						\cline{2-10}
						&\multirow{3}{*}{800}&961   &\multirow{3}{*}{1.78e-5}    &7 &4.27  &6.06e-4&6 &50.26 &6.06e-4\\
						&                    &3969  &                            &7 &17.58 &6.09e-4&6 &172.62&6.09e-4\\
						&                    &16129 &                            &7 &72.43 &6.09e-4&6 &649.47&6.09e-4\\		   
						
						\hline
						\multirow{9}{*}{1e-3}&\multirow{3}{*}{200}&961   &\multirow{3}{*}{2.85e-5}   &11&1.49  &1.38e-3&11&9.17 &1.38e-3\\
						&                    &3969  &    &11&6.41  &1.53e-3&11&26.04&1.53e-3\\
						&                    &16129 &   &11&26.94 &1.57e-3&11&95.74&1.57e-3\\
						\cline{2-10}
						&\multirow{3}{*}{400}&961   &\multirow{3}{*}{7.13e-6}   &12&3.18 &5.85e-4 &10 &23.97 &5.85e-4\\
						&                    &3969  &   &11&13.17&7.41e-4 &11 &90.67 &7.41e-4\\
						&                    &16129 &   &11&52.87&7.80e-4 &11 &325.89&7.80e-4\\
						\cline{2-10}
						&\multirow{3}{*}{800}&961   &\multirow{3}{*}{1.78e-6}   &12&6.75  &1.88e-4&10 &75.18  &1.88e-4\\
						&                    &3969  &   &11&26.09 &3.44e-4&11 &315.21 &3.44e-4\\
						&                    &16129 &  &11&107.39&3.83e-4&11 &1182.51&3.83e-4\\		                    
						\hline
						\multirow{9}{*}{1e-1}&\multirow{3}{*}{200}&961   &\multirow{3}{*}{2.85e-6}   &7 &1.01  &6.16e-4&7 &5.90 &6.16e-4\\
						&                    &3969  &    &7 &4.28  &1.24e-4&7 &16.91&1.24e-4\\
						&                    &16129 &   &7 &17.92 &1.20e-4&7 &61.84&1.20e-4\\
						\cline{2-10}
						&\multirow{3}{*}{400}&961   &\multirow{3}{*}{7.13e-7}   &8 &2.24  &6.43e-4&7 &16.89 &6.43e-4\\
						&                    &3969  &   &7 &8.72  &1.41e-4&7 &58.28 &1.41e-4\\
						&                    &16129 &  &7 &35.57 &6.07e-5&7 &208.84&6.07e-5\\
						\cline{2-10}
						&\multirow{3}{*}{800}&961   &\multirow{3}{*}{1.78e-7}   &8 &4.73  &6.57e-4&7 &52.66 &6.57e-4\\
						&                    &3969  &  &7 &17.40 &1.54e-4&7 &201.13&1.54e-4\\
						&                    &16129 &  &7 &72.01 &3.10e-5&7 &756.85&3.10e-5\\		                    
						\hline
						\multirow{9}{*}{1e1} &\multirow{3}{*}{200}&961   &\multirow{3}{*}{2.85e-7}   &4 &0.66  &6.82e-4&4 &3.44 &6.82e-4\\
						&                    &3969  &    &4 &2.70  &1.68e-4&4 &10.16&1.68e-4\\
						&                    &16129 &   &4 &11.18 &1.22e-4&4 &36.49&1.22e-4\\
						\cline{2-10}
						&\multirow{3}{*}{400}&961   &\multirow{3}{*}{7.13e-8}   &4 &1.33  &6.84e-4&4 &9.87  &6.84e-4\\
						&                    &3969  &   &4 &5.50  &1.70e-4&4 &33.78 &1.70e-4\\
						&                    &16129 &   &4 &22.33 &6.15e-5&4 &122.00&6.15e-5\\
						\cline{2-10}
						&\multirow{3}{*}{800}&961   &\multirow{3}{*}{1.78e-8}   &4 &3.66  &6.85e-4&4 &30.62 &6.85e-4\\
						&                    &3969  &   &4 &11.10 &1.71e-4&4 &116.08&1.71e-4\\
						&                    &16129 &  &4 &45.33 &4.25e-5&4 &436.58&4.25e-5\\		                    
						\hline
					\end{tabular}
				\end{center}
			\end{table}	
			
		}
	\end{example}
	\subsection{The Application of ${\bf P}_{\alpha}$ to A Locally Controlled Optimal Control Problem}
	In this subsection, we apply the proposed PinT preconditioner ${\bf P}_{\alpha}$ to the Schur complement system arising from a locally controlled optimal control problem, and present the related numerical results. The locally controlled optimal control problem is defined as follows
	\begin{equation}\label{optloccontrlprob}
		\min\limits_{y,u} L(y,u):=\frac{1}{2}||y-g||_{L^2(\Omega\times(0,T))}^2+\frac{\gamma}{2}||u||_{L^2(\Omega_0\times(0,T))}^2,
	\end{equation}
	subject to a linear wave equation with initial- and boundary-value  conditions 
	\begin{equation*}
		\begin{cases}
			y_{t}- Ly=f+\chi_{\Omega_0}u, \quad {\rm in~~}\Omega\times(0,T),\\
			y=0,\quad{\rm on~~}\partial\Omega\times(0,T),\\
			y(\cdot,0)=y_0,\quad \Omega.
		\end{cases}
	\end{equation*}
	The difference between the optimal control problem \eqref{optcontrlprob} and the locally controlled one \eqref{optloccontrlprob} is that  the regularization term $\frac{\gamma}{2}||u||_{L^2(\Omega_0\times(0,T))}^2$ is imposed on a subdomain $\Omega_0$ with $\Omega_0\subset\Omega$ and ${\rm Vol}(\Omega_0)>0$. Moreover, the control variable in PDE constraint of \eqref{optloccontrlprob} is weighted by the characteristic function $\chi_{\Omega_0}$. Here, $\chi_{\Omega_0}$ denotes the characteristic function of $\Omega_0$, which is defined as follows
	\begin{equation*}
		\chi_{\Omega_0}(x):=\begin{cases}
			1,\quad  x\in\Omega_0,\\
			0,\quad {\rm otherwise}.
		\end{cases}
	\end{equation*} 
	Similar to the discussion in Section \ref{introduction} for \eqref{optcontrlprob}, we can derive a reduced KKT system for \eqref{optloccontrlprob} as follows
	\begin{equation}\label{lockkt}
		\left[
		\begin{array}[c]{cc}
			I& L_2\\
			L_1&-\gamma^{-1}\chi_{\Omega_0}
		\end{array}
		\right]\left[\begin{array}[c]{c}
			y\\
			p
		\end{array}\right]=\left[\begin{array}[c]{c}
			g\\
			f
		\end{array}\right].
	\end{equation}
	Applying the same discretization scheme as that used for approximating \eqref{reducedkkt}, we obtain the following discrete KKT system for \eqref{lockkt}
	\begin{equation}\label{locdisckkt}
		\left[\begin{array}[c]{cc}
			\frac{\tau}{2} B_2\otimes I_J&B_1^T\otimes I_J+\frac{\tau}{2}B_2^T\otimes L_{h}\\
			B_1\otimes I_J+\frac{\tau}{2}B_2\otimes L_{h}&-\frac{\tau}{2\gamma}B_2^{\rm T}\otimes I_{\Omega_0}
		\end{array}\right]\left[\begin{array}[c]{c}
			y_{\tau,h}\\
			p_{\tau,h}
		\end{array}\right]=\left[\begin{array}[c]{c}
			g_{\tau,h}\\
			f_{\tau,h}
		\end{array}\right].
	\end{equation}
	Except for $I_{\Omega_0}$, other notations appearing in \eqref{locdisckkt} are already defined in \eqref{origdiscrdckktsys}. Here, $I_{\Omega_0}$ is a diagonal matrix whose $i$th diagonal element is defined as follows
	\begin{equation*}
		I_{\Omega_0}(i,i):=\begin{cases}
			1,\quad {\rm if~}i{\rm th~spatial~grid~point~locates~in~}\Omega_0,\\
			0,\quad {\rm if~}i{\rm th~spatial~grid~point~locates~in~}\Omega\setminus\Omega_0.
		\end{cases}
	\end{equation*}
	Applying the same block-diagonal scaling technique used in \eqref{scaledkktsys1} to \eqref{locdisckkt}, we find that solving the linear system \eqref{locdisckkt} reduces to solving the following Schur complement system
	\begin{equation}\label{locschurcmpsys}
		(\tau I_N\otimes I_{\Omega_0}+\eta GG^{\rm T})\tilde{v}=\tilde{b},
	\end{equation}
	for some given vector $\tilde{b}\in\mathbb{R}^{NJ\times 1}$. Clearly, \eqref{locschurcmpsys} is a real symmetric positive definite system, as $G$ (see the definition in \eqref{gmatdef}) is non-singular. Hence, both PCG-${\bf P}_{\alpha}$ and PCG-${\bf P}$ are applicable to solving \eqref{locschurcmpsys}.
	\begin{example}\label{loccontrolexpl}
		{\rm 
			In this example, we consider the locally controlled optimal control problem \eqref{optloccontrlprob} with
			\begin{align*}
				&a(x_1,x_2)\equiv 1,\quad\Omega=(0,1)^2,\quad T=1,\quad f(x_1,x_2,t)=(2\pi^2-1)\sin(\pi x_2)\sin(\pi x_1)\exp(-t),\\
				&g(x_1,x_2,t)=\sin(\pi x_2)\sin(\pi x_1)\exp(-t),\quad\Omega_0:=(0,1)^2\setminus(0,0.5)^2,
			\end{align*}
			the analytical solution of which is given by
			\begin{align*}
				y(x_1,x_2,t)=\sin(\pi x_1)\sin(\pi x_2)\exp(-t),\quad u\equiv 0.
			\end{align*}
			We apply PCG-${\bf P}_{\alpha}$ and PCG-${\bf P}$ to solving Example \ref{loccontrolexpl}, the numerical results of which are presented in Table \ref{loccontrolexpl}. From Table \ref{loccontrolexpl}, we see that the CPU cost pf PCG-${\bf P}_{\alpha}$ is much smaller than that of PCG-${\bf P}$ while the iteration numbers and errors of the two solvers are roughly the same. This again demonstrates the superiority of the proposed PinT preconditioner over the non-PinT preconditioner ${\bf P}$ in terms of computational time. Another observation from Table \ref{locexpltab} is that the iteration number of the both solvers decreases as $\gamma$ increases. The reason may be explained as follows. Both ${\bf P}_{\alpha}$ and ${\bf P}$ are originally designed for approximating the matrix $(\tau I_N\otimes I_J+\eta GG^{\rm T})$ while in Example \ref{loccontrolexpl}, ${\bf P}_{\alpha}$ and ${\bf P}$  are used as preconditioners of $(\tau I_N\otimes I_{\Omega_0}+\eta GG^{\rm T})$. Recall that $\eta=\gamma/\tau$. When $\gamma$ increases, $\eta GG^{\rm T}$ becomes dominant and $(\tau I_N\otimes I_{\Omega_0}+\eta GG^{\rm T})$ behaves more like $(\tau I_N\otimes I_{x}+\eta GG^{\rm T})$ in spectral sense. Therefore, when $\gamma$ is not small, ${\bf P}_{\alpha}$ and ${\bf P}$ approximate $(\tau I_N\otimes I_{\Omega_0}+\eta GG^{\rm T})$ well in spectral sense and both solvers converges quickly in such case. Hence, ${\bf P}_{\alpha}$ is an efficient preconditioner for locally controlled optimal control problem when $\gamma$ is not small.
			
			\newpage
			\begin{table}[H]
				\begin{center}
					\caption{Performance of the PCG-${\bf P}_{\alpha}$ and PCG-${\bf P}$ on Example \ref{loccontrolexpl}.
					}\label{locexpltab}
					\setlength{\tabcolsep}{0.7em}
					\renewcommand{\arraystretch}{0.9}
					\begin{tabular}[c]{cccc|ccc|ccc}
						\hline
						\multirow{2}{*}{$\gamma$} &\multirow{2}{*}{$N$} &\multirow{2}{*}{$J$}&\multirow{2}{*}{$\alpha$}&\multicolumn{3}{c|}{PCG-${\bf P}_{\alpha}$}&\multicolumn{3}{c}{PCG-${\bf P}$}\\
						\cline{5-10}
						&&&&Iter & CPU(s)& ${\rm E}_{N,J}$&Iter & CPU(s)& ${\rm E}_{N,J}$\\
						\hline
						\multirow{9}{*}{1e-4}&\multirow{3}{*}{100}&961   &\multirow{3}{*}{3.61e-4}    &24 &1.60 &4.61e-3 &23 &8.44 &4.61e-3\\
						&                    &3969  &    &23 &7.15 &4.63e-3 &23 &20.32&4.63e-3\\
						&                    &16129 &   &23 &29.50&4.63e-3 &23 &64.71&4.63e-3\\
						\cline{2-10}
						&\multirow{3}{*}{200}&961   &\multirow{3}{*}{9.02e-5}   &25 &3.59 &2.29e-3 &23 &20.11 &2.29e-3\\
						&                    &3969  &    &24 &15.23&2.30e-3 &23 &59.25 &2.30e-3\\
						&                    &16129 &   &24 &60.85&2.31e-3 &23 &206.51&2.31e-3\\
						\cline{2-10}
						&\multirow{3}{*}{400}&961   &\multirow{3}{*}{2.26e-5}   &25 &7.54  &1.14e-3 &23 &62.20 &1.14e-3\\
						&                    &3969  &   &25 &31.23 &1.15e-3 &23 &199.52&1.15e-3\\
						&                    &16129 &   &25 &128.73&1.15e-3 &23 &716.28&1.15e-3\\	
						\hline
						\multirow{9}{*}{1e-3}&\multirow{3}{*}{100}&961   &\multirow{3}{*}{1.14e-4}    &15 &1.14  &2.72e-3 &14 &4.78  &2.72e-3\\
						&                    &3969  &    &15 &4.77  &2.90e-3 &14 &12.06 &2.90e-3\\
						&                    &16129 &   &15 &19.87 &2.95e-3 &14 &39.27 &2.95e-3\\
						\cline{2-10}
						&\multirow{3}{*}{200}&961   &\multirow{3}{*}{2.85e-5}   &15 &2.24  &1.25e-3 &14 &12.22 &1.24e-3\\
						&                    &3969  &    &15 &9.75  &1.42e-3 &14 &36.04 &1.42e-3\\
						&                    &16129 &   &15 &39.32 &1.47e-3 &14 &124.66&1.47e-3\\
						\cline{2-10}
						&\multirow{3}{*}{400}&961   &\multirow{3}{*}{7.13e-6}   &15 &4.66  &5.23e-4 &14 &39.80 &5.23e-4\\
						&                    &3969  &   &15 &19.60 &6.83e-4 &14 &127.06&6.83e-4\\
						&                    &16129 &   &15 &79.89 &7.27e-4 &14 &426.76&7.27e-4\\		   
						
						\hline
						\multirow{9}{*}{1e-2}&\multirow{3}{*}{100}&961   &\multirow{3}{*}{3.61e-5}    &11 &0.77  &2.40e-4 &11 &3.75 &2.40e-4\\
						&                    &3969  &   &11 &3.60  &5.59e-4 &11 &9.50 &5.59e-4\\
						&                    &16129 &   &11 &15.00 &6.68e-4 &11 &30.80&6.68e-4\\
						\cline{2-10}
						&\multirow{3}{*}{200}&961   &\multirow{3}{*}{9.02e-6}   &11 &1.65  &3.35e-4 &11 &9.59 &3.35e-4\\
						&                    &3969  &    &11 &7.38  &2.09e-4 &11 &28.26&2.09e-4\\
						&                    &16129 &   &11 &29.84 &3.14e-4 &11 &96.93&3.14e-4\\
						\cline{2-10}
						&\multirow{3}{*}{400}&961   &\multirow{3}{*}{2.26e-6}   &11 &3.47  &4.54e-4 &11 &28.83 &4.54e-4\\
						&                    &3969  &  &11 &14.69 &5.98e-5 &11 &91.70 &5.98e-5\\
						&                    &16129 &   &11 &59.93 &1.39e-4 &11 &330.21&1.39e-4\\		                    
						\hline
						\multirow{9}{*}{1e-1}&\multirow{3}{*}{100}&961   &\multirow{3}{*}{1.14e-5}    &7  &0.59  &5.93e-4 &8 &2.76 &5.93e-4\\
						&                    &3969  &    &7  &2.35  &2.43e-4 &8 &6.91 &2.43e-4\\
						&                    &16129 &   &7  &9.99  &2.40e-4 &8 &22.93&2.40e-4\\
						\cline{2-10}
						&\multirow{3}{*}{200}&961   &\multirow{3}{*}{2.85e-6}   &7 &1.11   &6.36e-4 &8 &6.82 &6.36e-4\\
						&                    &3969  &    &7 &4.88   &1.30e-4 &8 &20.69&1.30e-4\\
						&                    &16129 &   &7 &19.80  &1.21e-4 &8 &71.53&1.21e-4\\
						\cline{2-10}
						&\multirow{3}{*}{400}&961   &\multirow{3}{*}{7.13e-7}   &8 &2.61   &6.56e-4 &7 &18.27 &6.56e-4\\
						&                    &3969  &   &7 &9.93   &1.50e-4 &8 &67.01 &1.50e-4\\
						&                    &16129 &   &7 &43.19  &6.09e-5 &8 &241.76&6.09e-5\\		                    
						\hline
						\multirow{9}{*}{1}   &\multirow{3}{*}{100}&961   &\multirow{3}{*}{3.61e-6}    &5 &0.44   &6.67e-4 &6 &2.11 &6.67e-4\\
						&                    &3969  &    &5 &1.78   &2.45e-4 &6 &5.34 &2.45e-4\\
						&                    &16129 &   &5 &7.47   &2.42e-4 &6 &17.53&2.42e-4\\
						\cline{2-10}
						&\multirow{3}{*}{200}&961   &\multirow{3}{*}{9.02e-7}   &5 &0.84   &6.78e-4 &6 &5.27 &6.78e-4\\
						&                    &3969  &    &5 &3.66   &1.65e-4 &6 &15.63&1.65e-4\\
						&                    &16129 &   &5 &14.88  &1.22e-4 &6 &54.30&1.22e-4\\
						\cline{2-10}
						&\multirow{3}{*}{400}&961   &\multirow{3}{*}{2.26e-7}   &5 &1.78   &6.82e-4 &5 &13.38 &6.82e-4\\
						&                    &3969  &   &5 &7.34   &1.68e-4 &6 &50.78 &1.68e-4\\
						&                    &16129 &   &5 &30.16  &6.14e-5 &6 &182.79&6.14e-5\\		                    
						\hline
					\end{tabular}
				\end{center}
			\end{table}	
		}
	\end{example}		
	\section{Concluding Remarks}
	In this paper, a PinT preconditioner has been proposed for the Schur complement system arising from the parabolic PDE-constrained optimal control problem. Theoretically, we have shown that the spectrum of the preconditioned matrix is lower and upper bounded by positive constants independent of matrix size and  the regularization parameter, thanks to which PCG solver for the preconditioned system has been proven to have an optimal convergence rate. Numerical results reported have shown the efficiency of the proposed preconditioning technique and supported the theoretical results.
	
	\section*{Acknowledgements}
	Thanks  Prof. Zhang, Zhi-Min (Wayne State University) for his comments and suggestions, which led to significant improvement. The work of Xue-lei  Lin was partially supported by research grants: 2021M702281 from China Postdoctoral Science Foundation, 12301480 from NSFC,  HA45001143 from Harbin Institute of Technology, Shenzhen, HA11409084  from Shenzhen. The work of Shu-Lin Wu was  supported by Jilin Provincial Department of Science and Technology (No. YDZJ202201ZYTS593).
	
	\bibliographystyle{plainnat}
	\bibliography{myreference}
\end{document}